\documentclass[a4paper]{article}

\usepackage[sort&compress,numbers]{natbib}
\usepackage[utf8x]{inputenc}
\usepackage[german,english]{babel}

\usepackage{psfrag}
\usepackage{epsfig}
\usepackage{subfigure}
\usepackage{amssymb}
\usepackage{amsmath}
\usepackage{amsthm}
\usepackage{courier}
\usepackage{listings}
\usepackage{booktabs}
\usepackage[colorlinks=true]{hyperref}
\usepackage[b]{esvect}   
\usepackage{float}
\usepackage{multirow}
\newcommand{\R}{\mathbb{R}}

\newcommand{\Pol}{\mathbb{P}}

\newcommand{\norm}[1]{\left\|#1\right\|}
\newcommand{\half}{\frac{1}{2}}

\begin{document}

%opening
\title{Compression challenges in large scale PDE solvers}

% Author Orchid ID: enter ID or remove command
\newcommand{\orcidauthorA}{0000-0003-0287-2120} % Add \orcidA{} behind the author's name
\newcommand{\orcidauthorB}{0000-0002-1071-0044} % Add \orcidB{} behind the author's name

\author{Sebastian Götschel\footnotemark[1], Martin Weiser\footnote{Zuse Institute Berlin, Takustr. 7, 14195 Berlin, \{goetschel,weiser\}@zib.de}}

\maketitle
\pagestyle{myheadings}
\thispagestyle{plain}

\begin{abstract}
  Solvers for partial differential equations (PDE) are one of the cornerstones of computational science.
  For large problems, they involve huge amounts of data that needs to be stored and transmitted on
  all levels of the memory hierarchy. Often, bandwidth is the limiting factor due to relatively small arithmetic intensity, and increasingly so due to the growing disparity between computing power and bandwidth.
  Consequently, data compression techniques have been investigated and tailored towards the specific requirements of PDE solvers during the last decades.
  This paper surveys data compression challenges and discusses examples of corresponding solution approaches for PDE problems, covering all levels of the memory hierarchy from mass storage up to main memory. Exemplarily, we illustrate concepts at particular methods, and give references to alternatives.\\

  \noindent\textbf{Keywords:} partial differential equation, data compression, floating point compression, lossy compression \\
  \textbf{MSC 2010: 65-02, 65M60, 65N30, 68-02, 68P30}
\end{abstract}

\section{Introduction}

Partial differential equations (PDEs) describe many phenomena in the natural sciences. Due to the broad spectrum of applications in physics, chemistry, biology, medicine, engineering, and economics, ranging from quantum dynamics to cosmology, from cellular dynamics to surgery planning, and from solid mechanics to weather prediction, solving PDEs is one of the cornerstones of modern science and economics. The analytic solution of PDEs in the form of explicit expressions or series representations~is, however, only possible for the most simplistic cases. Numerical simulations using finite difference, finite element, and finite volume methods~\cite{Strikwerda2007,DeuflhardWeiser2012,ZienkiewiczTaylorZhu2005} approximate the solutions on spatio-temporal meshes, and are responsible for a significant part of the computational load of compute clusters and high performance computing facilities worldwide.

Achieving sufficient accuracy in large PDE systems and on complex geometries often requires huge numbers of spatial degrees of freedom (up to $10^9$) and many time steps. Thus, numerical simulation algorithms involve large amounts of data that need to be stored, at least temporarily, or transmitted to other compute nodes running in parallel. Therefore, large scale simulations face two main data-related challenges: communication bandwidth and storage capacity.

First, computing, measured in floating point operations per second (FLOPS), is faster than data transfer, measured in bytes per second. The ratio has been increasing for the last three decades, and continues to grow with each new CPU and GPU 
 generation~\cite{McCalpin2016}. Today, the performance of PDE solvers is mostly limited by communication bandwidth, with the CPU cores achieving only a tiny fraction of their peak performance. This concerns the CPU-memory communication, the so-called ``memory wall''~\cite{McCalpin1995,McKee2004,Alted2010}, as well as inter-node communication in large distributed systems~\cite{ReedDongarra2015}, and popularized the ``roofline model'' as a means to understand and interpret computer performance.

Second, storage capacity is usually a limited resource. Insufficient storage capacity can affect simulations in two different aspects. If needed for conducting the computation, it limits the size of problems that can be treated, and thus the accuracy of the results. Alternatively, data can spill over to the next larger and slower level of the memory hierarchy, with a corresponding impact on the simulation performance. If needed for storing the results, the capacity limits the number or resolution of simulation results that can be used for later interpretation, again affecting the accuracy of the conclusions drawn from the simulations. For both aspects, a variety of data compression methods have been proposed, both as pure software solutions and as improved hardware architecture.

The aspects of data compression that are specific for PDE solvers and the computed solutions are reviewed in the following Section~\ref{sec:aspects}, where we also classify compression methods proposed in the literature according to these aspects. After that, we discuss use cases of compression in PDE solvers along the memory hierarchy at prototypical levels of in-memory compression (Section~\ref{sec:memory}), inter-node communication (Section~\ref{sec:distributed}), and mass storage (Section~\ref{sec:storage}). For each use case, an example is presented in some detail along with the one or two compression methods applied, highlighting the different needs for compression on one hand and the different challenges and trade-offs encountered on the other~hand.

\section{Compression Aspects of PDE Solvers}\label{sec:aspects}
PDE solvers have a requirement profile for compression that differs in several aspects from other widespread compression demands like text, image, video, and audio compression. The data to be compressed consists mainly of raw floating point coefficient vectors in finite difference, finite element, and finite volume methods.

\subsection{High Entropy Data Necessitates Lossy Compression}
Usually, double precision is used for coefficient vectors in order to avoid excessive accumulation of rounding errors during computation, even if the accuracy offered by 53 mantissa bits is not required for representing the final result. Due to rounding errors, the less significant mantissa bits are essentially random, and incur a large entropy. Lossless compression methods are, therefore, not able to achieve substantial compression factors, i.e., ratios of original and compressed data sizes. Examples of lossless compression schemes tailored towards scientific floating point data are fpzip 
~\cite{LindstromIsenburg2006}, FPC~\cite{BurtscherRatanaworabhan2009}, SPDP\footnote{\url{https://userweb.cs.txstate.edu/~burtscher/research/SPDP/}}
~\cite{ClaggettAzimiBurtscher2018}, Blosc\footnote{\url{http://blosc.org/}}
~\cite{Alted2010}, and Adaptive-CoMPI~\cite{FilSinCarCalGar2011}.

In contrast, lossy compression allows much higher compression factors, but requires a~careful selection of compression error in order not to compromise the final result of the computation. In general, there is no need for the compression error to be much smaller than the discretization or truncation errors of the computation. Lossy compression schemes that have been proposed for scientific floating point data in different contexts include ISABELA\footnote{\url{http://freescience.org/cs/ISABELA/ISABELA.html}} (In-situ Sort-And-B-spline Error-bounded Lossy Abatement)
~\cite{LakshminarasimhanShahEthierKuChangKlaskyLathamRossSamatova2013}, SQE~\cite{IversonKamathKarypis2012}, zfp\footnote{\url{https://computation.llnl.gov/projects/floating-point-compression}}~\cite{Lindstrom2014,Lindstrom2017,DifFoxHitSanLin2018}, SZ\footnote{\url{https://collab.cels.anl.gov/display/ESR/SZ}}
 1.1~\cite{DiCap2016} and 1.4~\cite{TaoDiCheCap2017,LiaDiTaoLiLiGuoCheCap2018}, multilevel transform coding on unstructured grids (TCUG)~\cite{WsGoetschel2012,Goe2015}, adaptive thinning (AT)~\cite{DemaretDynFloaterIske2005,SolinAndrs2009} and adaptive coarsening (AC)~\mbox{\cite{ShafaatBaden2007,UnatHromadkaBaden2009}}, TuckerMPI~\mbox{\cite{AustinBallardKolda2016,BallardKlinvexKolda2019}}, TTHRESH~\cite{BallesterRipollLindstromPajarola2019}, MGARD~\cite{AinsworthTuglukWhitneyKlasky2019}, HexaShrink~\cite{PeyDuvPayBouChiSchAnt2019}, and hybrids of different methods~\cite{TaoDiLiaCheCap2019}.

\subsection{Data Layout Affects Compression Design Space}
One of the most important differences between PDE solvers from the compression point of view is the representation and organization of spatial data.

General-purpose floating point compression schemes essentially ignore the underlying structure and treat the values as an arbitrary sequence of values. Examples of this approach are ISABELA, SQE, and SZ-1.1. The advantage of direct applicability to any kind of discretization comes at the cost of moderate compression factors, since position-dependent correlations in the data cannot be directly~exploited.

Structured Cartesian grids are particularly simple and allow efficient random access by index computations, but are limited to quasi-uniform resolutions and relatively simple geometries parametrized over cuboids. Cartesian structured data is particularly convenient for compression, since it allows the use of the Lorenzo predictor (fpzip) or its higher order variants (SZ-1.4), regular coarsening~(AC), simple computation of multilevel decompositions (MGARD and HexaShrink), or exploiting tensor approaches such as factorized block transforms (zfp) or low-rank tensor approximations (TuckerMPI and TTHRESH).

If complex geometries need to be discretized or if highly local solution features need to be resolved, unstructured grids are used. Their drawback is that coefficients are stored in irregular patterns, and need to be located by lookup. Fewer methods are geared towards compression of data on unstructured grids. Examples include TCUG and AT. When storing values computed on unstructured grids, the grid connectivity needs to be stored as well. In most use cases, however, e.g., iterative solvers or time stepping, many coefficient vectors have to be compressed, but few grids do. Thus, the floating point data represents the bulk of the data to be compressed. Methods developed for compression of unstructured grid geometries can in some cases been used for storing solutions or coefficient data of PDE solvers, but are mostly tailored towards computer graphics needs. We refer to the survey~\cite{MagLavDupHud2015} for 3D mesh compression, as well as~\cite{GoetschelVTycowiczPolthierWeiser2015} for grid compression in PDE applications. In addition to geometry and connectivity, recent research focuses on the compression of attribute data, such as color values or texture, taking the geometry into account~\cite{NasBidPayMau2019}, and using progressive compression methods~\cite{CaiVidDupLav2016}, but without rigorous error control.

In contrast to the grid structure, the method of discretization (finite difference, finite element, or finite volume methods, or even spectral methods) is of minor importance for compression, but of course affects technical details.

PDE solvers can also benefit from compression of further and often intermediate data that arises during the solution process and may be of completely different structure. This includes in particular preconditioners (scalar quantization, mixed precision~\cite{AnztDongarraFlegarHighamQuintanaorti2019,SchneckWeiserWende2019}, hierarchical matrices~\cite{Hackbusch1999}), discretized integral operators in boundary element methods (wavelets~\cite{DahmenHarbrechtSchneider2006}), and boundary corrections in domain decomposition methods.

\subsection{Error Metrics and Error Propagation Affect Compression Accuracy Needs}
The notion of ``compression error'' is not well-defined, but needs to be specified in view of a~particular application. The impact of compression on data analysis has been studied empirically~\cite{LuLiuHeLuoSuchytaChoiPodhorszkiKlaskyWolfLiuQiao2018,AinsworthTuglukWhitneyKlasky2019} and statistically~\cite{PoppickNardiFeldmanBakerHammerling2018,BakHamMicXuEtAl2016}. Ideally, the compression scheme is tailored towards the desired error metric. In lack of knowledge about the application's needs, general error metrics such as pointwise error (maximum or $L^\infty$ norm) and mean squared error (MSE, $L^2$ norm) are ubiquitous, and generally used for compressor design. Integer Sobolev semi-norms and Hausdorff distances of level sets have been considered by Hoang et al.~\cite{HoangKlacanskyBhatiaBremerLindstromPascucci2019} for sorting level and bitplane contributions in wavelet compression. Broken Sobolev norms have been used by Whitney~\cite{Whitney2018} for multilevel decimation.

If intermediate values are compressed for storage, in addition to the final simulation results compressed for analysis and archival, errors from lossy compressions are propagated through the following computations. Their impact on the final result depends very much on the type of equation that is solved and on the position in the solution algorithm where the compression errors enter. For example, inexact iterates in the iterative solution of equation systems will be corrected in later iterations, but may lead to an increase in iteration count. The impact of initial value or source term errors on the solution of parabolic equations is described by negative Sobolev norms, since high-frequency components are damped out quickly. In contrast, errors affecting the state of an explicit time stepping method for hyperbolic equations will be propagated up to the final result. Such analytical considerations are, however, qualitative, and do not allow designing quantization tolerances for meeting a quantitative accuracy requirement.

Fortunately, a posteriori error estimates are often available and can be used for controlling compression errors as well as for rate-distortion optimization. As an example, error estimators have been used in~\cite{GoetschelWs2015} for adaptive selection of state compression in the adjoint computation of gradients for optimal control.

\subsection{Compression Speed and Complexity Follow the Memory Hierarchy}
Compression plays different roles on all levels of the memory hierarchy, depending on application, problem size, and computer architecture. 
Due to the speed of computation growing faster than memory, interconnect, and storage bandwidth, a multi-level memory hierarchy has developed, ranging from several memory cache levels over main memory, nonvolatile memory (NVRAM) and solid state disk burst buffers down to large storage systems~\cite{JacobNgWang2010}. Lower levels exhibit larger capacity, but less bandwidth than higher levels. The larger the data to be accessed, the deeper in the memory hierarchy it needs to be stored, and the slower is the access. Here, data compression can help to reduce the time to access the data and to exploit the available capacity on each level better. While this can, in~principle, be considered and tuned for any of the many levels of current memory hierarchies, we limit the discussion to three prototypical levels: main memory, interconnects, and storage systems.

Compression of in-memory data aims at avoiding the ``memory wall'' and reducing the run time of the simulation (see Section~\ref{sec:memory}). The available bandwidth is quite high, even if not sufficient for saturating the computed units. In order to observe an overall speedup, the overhead of compression and decompression must be very small, such that only rather simple compression schemes working on small chunks of data can be employed.

In distributed systems, compression of inter-node communication can be employed to mitigate the impact of limited network bandwidth on the run time of simulations (see Section~\ref{sec:distributed}). The bandwidth of communication links is about an order of magnitude below the memory bandwidth, and the messages exchanged are significantly larger than the cache lines fetched from memory, such that more sophisticated compression algorithms can be used. 

Mass storage comes into play when computed solutions need to be stored for archiving or later analysis. Here, data size reduction is usually of primal interest, such that complex compression algorithms exploiting correlations, both local and global, in large data sets can be employed (see ~Section~\ref{sec:storage}). For an evaluation of compression properties on several real-world data sets we refer to~\cite{LuLiuHeLuoSuchytaChoiPodhorszkiKlaskyWolfLiuQiao2018}. Due to the small available bandwidth and the correspondingly long time for reading or writing uncompressed data, the execution time even of complex compression algorithms can be compensated when storing only smaller compressed data sets. This aspect is relevant for the performance of out-of-core algorithms for very large problems. If compression data can be kept completely in memory, out-of-core algorithms can even be turned to in-core algorithms. A recent survey of use cases for reducing or avoiding the I/O bandwidth and capacity requirements in high performance computing, including results using mostly SZ and zfp, is given by Cappello et al.~\cite{CappelloDiLiLiangGokTaoYoonWuAlexeevChong2019}.

\section{In-Memory Compression} \label{sec:memory}

\begin{figure}[H]
\centering
 \includegraphics[width=0.7\textwidth]{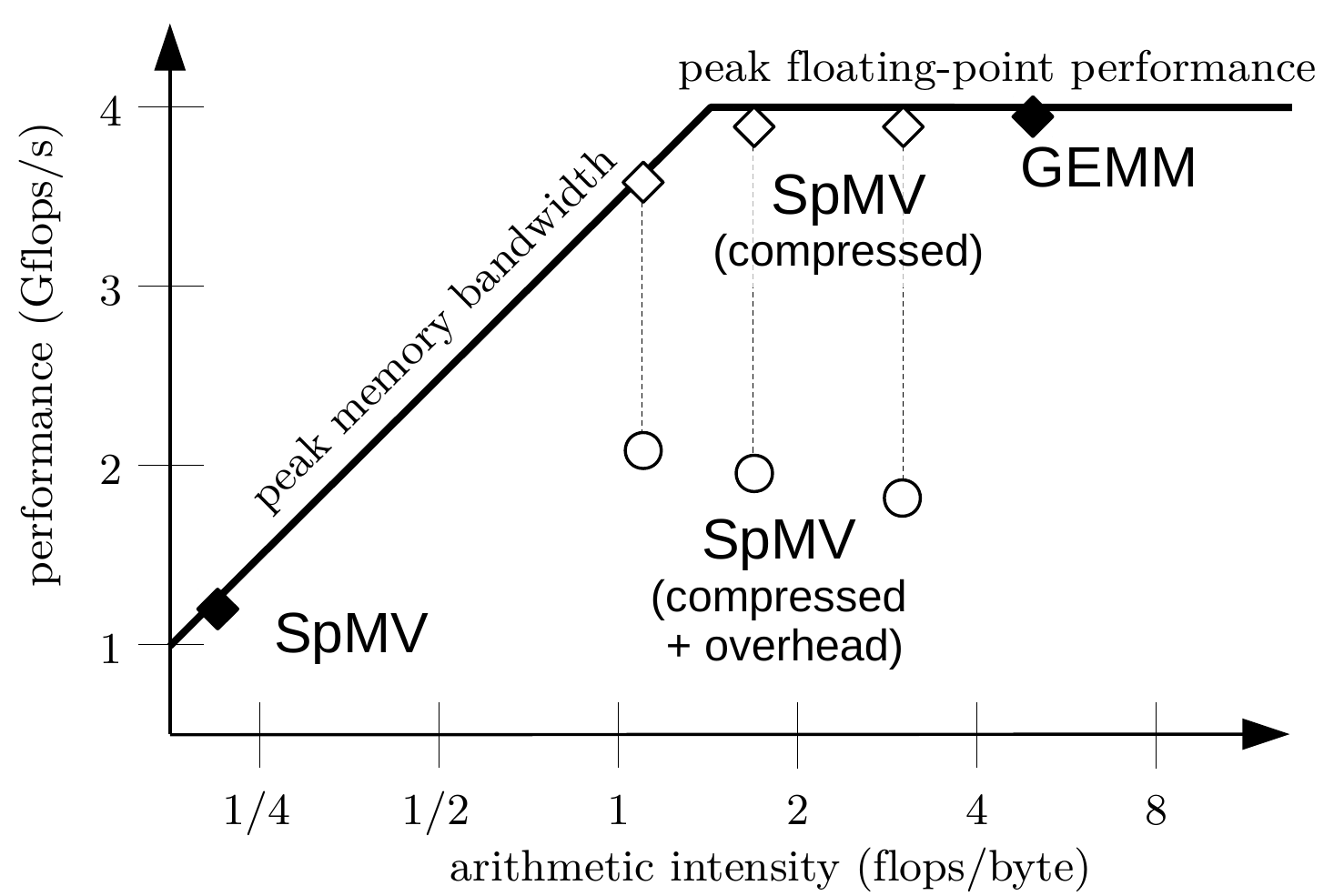}
 \caption{Naive roofline model showing achievable performance vs. the arithmetic intensity. Some computations, e.g., dense matrix-matrix multiplication (GEMM), perform many flops per byte fetched or written to memory, such that their execution speed is bounded by the peak floating point~performance. Others, such as sparse matrix-vector multiplication (SpMV), require many bytes to be fetched from memory for each flop, and are therefore memory-bound (filled diamonds).
 Data~compression methods for memory-bound computations can reduce the amount of data to be read or written, and therefore increase the arithmetic intensity. Different compression schemes can achieve different compression factors (empty diamonds on top) and thus different arithmetic intensities. The~computational overhead of compression and decompression can, however, reduce the performance gain (empty circles bottom), depending on the complexity of the compression method used. }
 \label{fig:roofline}
\end{figure}

The arithmetic density of a numerical algorithm, i.e., the number of floating point operations performed per byte that is read from or written to memory, is one of the most important properties that determines the actual performance. With respect to that quantity, the performance can be described by the roofline model~\cite{WilliamsWatermanPatterson2009}. It includes two main bounds, the peak performance, and the peak memory bandwidth; see Figure~\ref{fig:roofline}. Most PDE solvers are memory bound, in particular finite element methods working on unstructured grids and making heavy use of sparse linear algebra, but also stencil-based finite difference schemes in explicit time stepping codes.

Performance improvements can be obtained by increasing the memory bandwidth, e.g., exploiting NUMA architectures or using data layouts favoring contiguous access patterns, or by reducing the amount of data read from and written to the memory, increasing the arithmetic intensity. Besides larger caches, data compression is an effective means to reduce the memory traffic. Due to the reduced total data size, it can also improve cache hit rates or postpone the need for paging or out-of-core algorithms for larger problem instances.

While compression can reduce the memory traffic such that the algorithm becomes compute-bound, the overhead of compression and decompression realized in the software reduces the budget available for payload flops. This is illustrated in Figure~\ref{fig:roofline}. Sparse matrix-vector products are usually memory bound with an arithmetic intensity of less than 0.25 flops/byte. Data compression increases the arithmetic intensity and moves the computation towards the peak floating point roofline. For illustration, compression factors between 4 and 16 are shown, representing different compression schemes.  The (de)compression overhead of one to three flops per payload flop reduces the performance delivered to the original computation. Depending on the complexity of the compression scheme, and hence its computational overhead, the resulting performance can even be worse than before. This implies that to overcome the memory wall, only very fast, and therefore rather simple, compression schemes are beneficial.

Dedicated hardware support can raise the complexity barrier for compression algorithms, and several such approaches have been proposed. While completely transparent approaches~\cite{PekhimnkoSeshadriKimXinMutluGibbonsKozuchMowry2013,ShafieeTaassoriBalasubramonianDavis2014,YoungNairQureshi2017} must rely on lossless compression and will therefore achieve only small compression factors on floating point data, intrusive approaches can benefit from application-driven accuracy of floating point representations~\cite{JainHillLinKhanHaqueLaurenzanoMahlkeTangMars2016}. They need, however, compiler support and extended instruction set architectures, and cannot be realized on commodity systems. For a survey on hardware architecture aspects for compression we refer to~\cite{MittalVetter2015}.

An important issue is the transport of compression errors through the actual computation and the influence on the results. While it is hard to envisage quantization being guided by a posteriori error estimates due to their computational overhead, careful analysis can sometimes provide quantitative worst-case bounds. Such an example, an iterative solver, is presented in the following section. Iterative solvers are particularly well-suited candidates for in-memory compression for two reasons. First, they often form the inner loops of PDE solvers and therefore cause the most memory traffic that can benefit from compression. Second, many can tolerate a considerable amount of relative error while still converging to the correct result. Thus, the compression error trade-off is not between compression factor and accuracy, but between compression factor and iteration count.

\subsection{Scalar Quantization for Overlapping Schwarz Smooth\-ers}
One particularly simple method of data compression is a simple truncation of mantissa and exponent bits, i.e., using IEEE~754 single precision (4 bytes) instead of double precision (8 bytes) representations~\cite{IEEE754}, or even the half precision format (2 bytes) popularized by recent machine learning applications. Conversion between the different formats is done in hardware on current CPUs and integrated into load/store operations, such that the compression overhead is minimal. Consequently, using mixed precision arithmetics has been considered for a long time, in particular in dense linear algebra~\cite{BaboulinEtAl2009} and iterative solvers~\cite{AnztLuszczekDongarraHeuveline2012,Grout2015}. Depending on the algorithm's position in the roofline model, either the reduced memory traffic (Basic Linear Algebra Subroutines (BLAS) level 1/2, vector-vector and matrix-vector operations) or the faster execution of lower precision floating point operations (BLAS level 3, matrix-matrix operations) is made use of.

An important building block of solvers for elliptic PDEs of the type
\begin{align*}
 -\mathop\mathrm{div}(\sigma\nabla u) &= f && \text{in $\Omega$} \\
 n^T \sigma \nabla u + \alpha u &= \beta && \text{on $\partial\Omega$}
\end{align*}
is the iterative solution of the sparse, positive definite, and ill-conditioned linear equation systems $Ax=b$ arising from finite element discretizations. For this task, usually preconditioned conjugate gradient (CG) methods are employed, often combining a multilevel preconditioner with a Jacobi smoother~\cite{DeuflhardWeiser2012}. For higher order finite elements, with polynomial ansatz order $p>2$, the effectivity of the Jacobi smoother quickly decreases, leading to slow convergence. Then it needs to be replaced by an overlapping block Jacobi smoother $B$, with the blocks consisting of all degrees of freedom associated with cells around a grid vertex. Application of this smoother then involves a large number of essentially dense matrix-vector multiplications of moderate size:
\begin{equation}
 B^{-1} = \sum_{\xi\in\mathcal{N}} P_\xi A_\xi^{-1} P_\xi^T.
\end{equation}
Here, $\mathcal{N}$ is the set of grid vertices, $A_\xi$ is the symmetric submatrix of $A$ corresponding to the vertex~$\xi$, and $P_\xi$ distributes the subvector entries into the global vector. Application of this smoother dominates the solver run time, and is strictly memory bound due to the large number of dense matrix-vector~multiplications.

Compressed storage of $A_\xi^{-1}$ as $\tilde A_\xi^{-1}$ using low precision representation of its entries has been investigated in~\cite{SchneckWeiserWende2019}. A detailed analysis reveals that the impact on the preconditioner effectivity and hence the CG convergence is determined by $\|A^{-1}-\tilde A^{-1}\|_2$. This suggests that a uniform quantization of submatrix entries should be preferable in view of rate-distortion optimization. Accordingly, fixed point representations have been considered as an alternative to low precision floating point representations. Moreover, the matrix entries exhibit a certain degree of correlation, which can be exploited by dividing $A_\xi^{-1}$ into square blocks to be stored independently. A uniform quantization of their entries $c$ within the block entries' range $[c_{\min},c_{\max}]$ as $l=q(c)$ and corresponding dequantizing  $\tilde c=q^+(l)$ is given by
\begin{equation}
 q(c) = \begin{cases} \lfloor \frac{c-c_{\rm min}}{\Delta}\rfloor, & c < c_{\rm max} \\
                         2^k -1, & c=c_{\rm max},
           \end{cases} \qquad
 q^+(l) =c_{\rm min} + \Delta \Big(l+\frac{1}{2}\Big),
\end{equation}
with step size $\Delta = 2^{-k}(c_{\rm max}-c_{\rm min})$, providing minimal entry-wise quantization error $\Delta/2$ for the given bit budget.

\begin{figure}[H]
\centering
 \includegraphics[width=\textwidth]{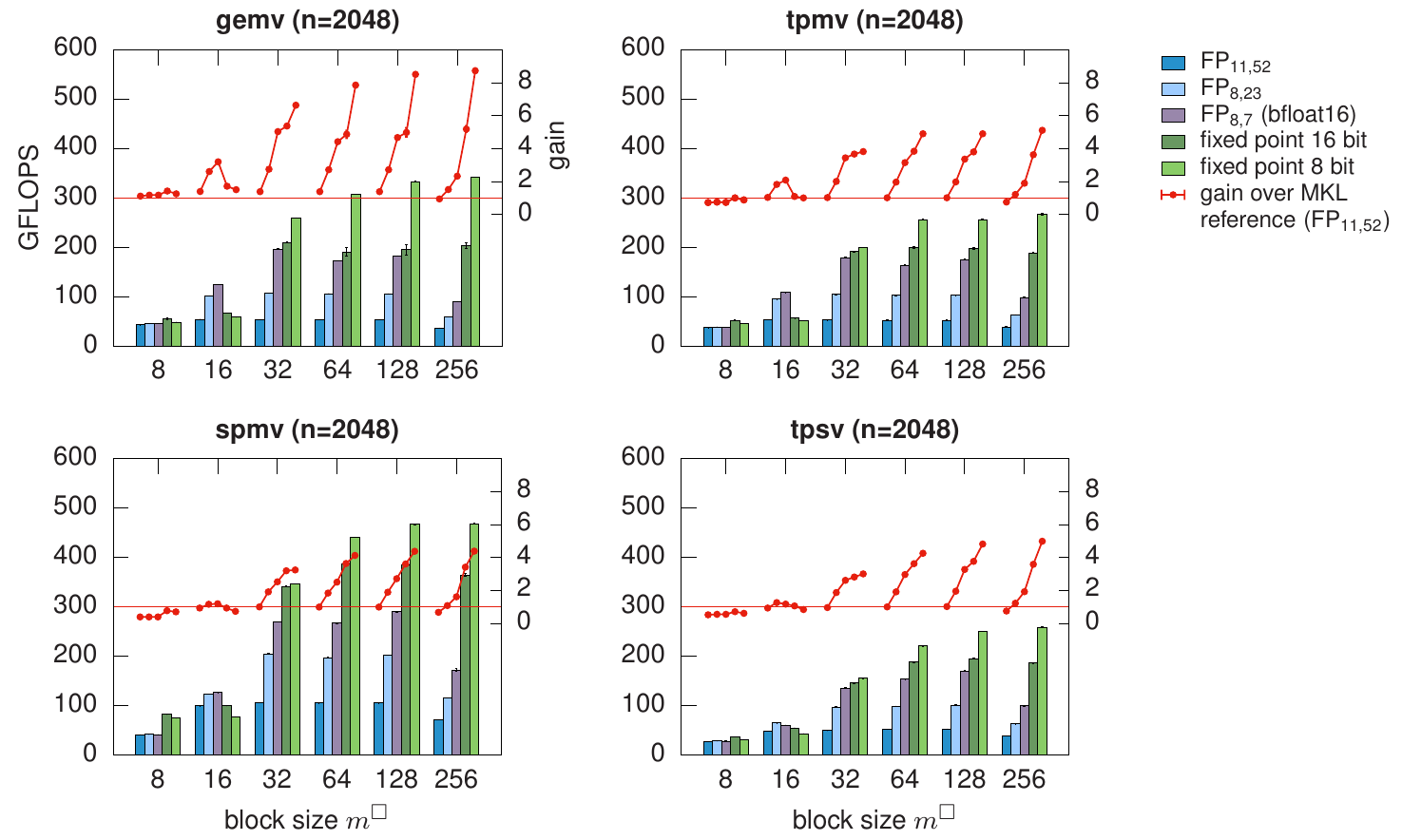}
 \caption{Run times of BLAS
 	 level 2 operations on $2048\times 2048$-matrices for overlapping Schwarz smoothers with mixed precision. Depending on the access patterns, a speedup over the Intel Math Kernel Library (MKL) almost on par with the compression factor can be achieved~\cite{SchneckWeiserWende2019}.}
 \label{fig:blas2}
\end{figure}

Decompression can then be performed inline during application of the preconditioner, i.e., during the matrix-vector products. The computational overhead is sufficiently small as long as conversions between arithmetic data types are performed in hardware, which restricts the possible compression factors to $\{2,4,8\}$, for which the speedup reaches almost the compression factor; see~Figure~\ref{fig:blas2}. With~direct hardware support for finer granularity of arithmetic data types to be stored~\cite{JainHillLinKhanHaqueLaurenzanoMahlkeTangMars2016}, an even better fit of compression errors to the desired accuracy could be achieved with low overhead.

The theoretical error estimates together with typical condition numbers of local matrices $A_\xi$ suggest that using a 16 bit fixed-point representation should increase the number of CG iterations by not more than 10\% due to preconditioner degradation, up to an ansatz order $p=5$. In fact, this is observed in numerical experiments, leading to a speedup of preconditioner application by a factor of up to four. With moderate ansatz order $p\le 6$, even 8 bit fixed point representations can be used, achieving a speedup of up to six~\cite{SchneckWeiserWende2019}.

Similar results have been obtained for non-overlapping block-Jacobi preconditioners for Krylov methods applied to general sparse systems~\cite{AnztDongarraFlegarHighamQuintanaorti2019} and for substructuring domain decomposition methods~\cite{GiraudHaidarWatson2006}.

We would like to stress that the bandwidth reduction is the driving motivation rather than the possible speedup due to faster single precision arithmetics, in contrast to BLAS level 3 algorithms. Not only is the preconditioner application memory bound such that the data size is the bottleneck, but there is also a compelling mathematical reason for performing the actual computations in high precision arithmetics: Storing the inverted submatrices $A_\xi^{-1}$ in low precision results in a valid, though less effective, symmetric positive definite preconditioner as long as the individual submatrices remain positive definite. Performing the dense matrix-vector products in low precision, however, will destroy the preconditioner's symmetry sufficiently to affect CG convergence, and therefore leave the well-understood theory of subspace correction methods.

\subsection{Fixed-Rate Transform Coding}
To achieve compression factors that are higher than possible with reduced precision storage, more sophisticated and computationally more expensive approaches are required. Competing aims are high compression factors, low computational overhead, as well as transparent and random access. One~particularly advanced approach is transform coding on Cartesian grids~\cite{Lindstrom2014}, which is the core of~zfp. Such structured grids, though restrictive, are used in those areas of scientific computing where no complex geometries have to be respected and the limited locality of solution features does not reward the overhead of local mesh refinement.

The straightforward memory layout of the data allows considering tensor blocks of values that are compressed jointly. For 3D grids, $4\times 4\times 4$ blocks appear to be a reasonable compromise between~locality, which is %Please confirm meaning is retained. -- OK
necessary for random access, and compression factors due to exploitation of spatial correlation. These blocks are transformed by an orthogonal transform. While well-known block transforms such as the discrete cosine transform~\cite{AhmNatRao1974} can be used, a special transform with slightly higher decorrelation efficiency has been developed in~\cite{Lindstrom2014}. Such orthogonal transforms can be applied efficiently by exploiting separability and lifting scheme for factorization, i.e., applying 1D transforms in each dimension, and realizing these 1D transforms by a sequence of cheap in-place modifications. This results in roughly 11 flops per coefficient. The transform coefficients are then coded bitplane by bitplane using group testing, similar to set partitioning in hierarchical trees~\cite{SaidPearlman1996}. This embedded coding schemes allows decoding data at variable bit rate, despite the fixed-rate compression enforced by random access ability.

Despite a judicious choice of algorithm parameters, which allow an efficient integer implementation of the transform using mainly bit shifts and additions, the compression and decompression are heavily compute-bound. The effective single-core throughput, depending on the compression factor, is reported to lie around 400 MB/s, which is about a factor of ten below contemporary memory bandwidth. Here, dedicated hardware support for transform coding during read and write operations would be beneficial. In order to be flexible enough to support different coding schemes, transforms, and data sizes, however, this would need to be configurable.

In conclusion, simple and less effective compression schemes such as mixed precision approaches appear to be today's choice for addressing the memory wall in PDE computations. Complex and more effective schemes are currently of interest mainly to fit larger problems into a given memory budget. This is, however, likely to change in the future: As the hardware continues to trend to more cores per CPU socket, and thus the gap between computing performance and memory bandwidth widens, higher complexity of in-memory compression will pay off. But even then, adaptive quantization based on a posteriori error estimates for the impact of compression error will probably be out of reach.

\section{Communication in Distributed Systems} \label{sec:distributed}

The second important setting in which data compression plays an increasingly important role in PDE solvers is communication in distributed systems. The ubiquitous approach for distribution is to partition the computational domain into several subdomains, which are then distributed to the different compute nodes. Due to the locality of interactions in PDEs, communication happens at the boundary shared by adjacent subdomains. A prime example are domain decomposition solvers for elliptic problems~\cite{ToselliWidlund2005}.

The inter-node bandwidth in such systems ranges from around 5 GB/s per link with high-performance interconnects such as InfiniBand down to shared 100 MB/s in clusters made of commodity hardware such as gigabit ethernet. This is about one to two orders of magnitude below the memory bandwidth. Consequently, distributed PDE solvers need to have a much higher arithmetic density with respect to inter-node communication than with respect to memory access. In~domain decomposition methods, the volume of subdomains in $\R^d$ with diameter $h$ scales with $h^d$, as does the computational work per subdomain. The surface and hence communication, however, scale only with $h^{d-1}$, such that high arithmetic intensity can be achieved by using sufficiently large subdomains---which impedes on weak scaling and limits the possible parallelism. Consequently, communication can become a severe bottleneck.

Data compression has been proposed for increasing the effective bandwidth. Burtscher and Ratanaworabhan~\cite{BurtscherRatanaworabhan2009} consider lossless compression of floating point data streams, focusing on high throughput due to low computational overhead. Combining two predictors based on lookup tables trained online from already seen data results in compression factors on par with other lossless floating point compression schemes and general-purpose codes like GZIP, at a vastly higher throughput. Being lossless and not exploiting the spatial correlation of PDE solution values limits the compression factor, however, to values between 1.3 and 2.0, depending on the size of the lookup tables. \mbox{Filgueira et al.~\cite{FilSinCarCalGar2011}} present a transparent compression layer for MPI communication, choosing adaptively between different lossless compression schemes. Again, with low redundancy of floating point data, as is characteristic for PDE coefficients, compression factors below two are achieved.

Higher compression factors can be achieved with lossy compression. As in in-memory compression, analytical a priori error estimates can provide valuable guidance on the selection of quantization tolerances. Here, however, the inter-node communication bandwidth is relatively small, such that the computational cost of a posteriori error estimators might be compensated by the additional compression opportunities they can reveal---an interesting topic for future research.

\subsection{Inexact Parallel-in-Time Integrators} \label{sec:PinT}

An example of communication in distributed systems is the propagation of initial values in parallel-in-time integrators for initial value problems $\dot u = f(u)$, $u(t_0) = a$, in particular of hybrid parareal type~\cite{Gander2015}. Here, the initial value problem is interpreted as large  equation system
\begin{equation}
F(U) = \begin{bmatrix}
 a &- u^0(t_0)  \\
 & \dot u^0 - f(u^0)  \\
 & u^0(t^1) &-  u^1(t^1 )\\
 && \dot u^1 - f(u^1) \\
 && u_1(t^2) &-  u^2(t^2 )\\
 &&&\ddots
 \end{bmatrix} = 0
\end{equation}
for a set $U=(u^0,\dots,u^N)$ of subtrajectories $u^n \in C^1([t^n,t^{n+1}])$ on a time grid $t^0,t^1, \dots, t^{N+1}$. Instead of the inherently sequential triangular solve, i.e., time stepping, the system is  solved by a stationary iterative method with an approximate solver $S$:
\begin{equation}
 U_{j+1} = U_{j} + S(F(U_j)).
\end{equation}
The advantage is that a large part of the approximate solver $S$ can be parallelized, by letting $S^{-1} = \mathcal{F}^{-1} + \mathcal{G}^{-1}$, where $\mathcal{G}$ is an approximation of the derivative $-F'$ on a spatial and/or temporal coarse grid and provides the global transport of information, while $\mathcal{F}$ is a block-diagonal approximation of $-F'$ on the fine grid and cares for the local reduction of fine grid residuals. The bulk of the work is done in applying the fine grid operator $\mathcal{F}^{-1}$, where all blocks can be treated in parallel. Only the coarse grid solution operator $\mathcal{G}^{-1}$ needs to be applied sequentially. If the parallelized application of $S$ is significantly faster than computing a single subtrajectory up to fine grid discretization accuracy, reasonable parallel efficiencies above 0.5 can be achieved~\cite{EmmettMinion2012}.

For a fast convergence, however, the terminal values $u^n(t^{n+1})$ have to be propagated sequentially as initial values of $u^{n+1}(t^{n+1})$ over all subintervals during each application of the approximate solver~$S$. Thus, communication time can significantly affect the overall solution time~\cite{FischerGoetschelWs2018}. Compressed communication can therefore improve the time per iteration, but may also impede on the convergence speed and increase the number of iterations. A judicious choice of compression factor and distortion must rely on error estimates and run time models.

The worst-case error analysis presented in~\cite{FischerGoetschelWs2018} provides a bound of the type
\begin{equation}
 \|U_j - U_*\| \le c_j^n \, \left(\frac{1+\Delta_C}{1-\Delta_C/\rho}\right)^{n+2},
\end{equation}
depending on the relative compression error $\Delta_C$, the local contraction rate $\rho$ of $S$, and factors $c_j^n$ independent of communication. This can be used to compute an upper bound on the number $J(\Delta_C)$ of iterations in dependence of the compression error. The run time of the whole computation is $T_{\rm par} = N(t_G+t_C(\Delta_C)) + J(\Delta_C)(t_G+t_F)$, where $t_G$ is the time required for the sequential part of $S$, $t_F$ for the parallel part, and the function $t_C$ is the communication time depending on the compression error, including compression time. Minimizing $T_{\rm par}$ can be used to optimize for $\Delta_C$, as long as the relation between $\Delta_C$ and $t_C$ is known. For finite element coefficients, most schemes will lead to $t_C \approx  -c \log \Delta_C$, i.e., a communication time proportional to the bits spent per coefficient, with the proportionality factor $c$ depending on bandwidth, problem size, and efficiency of the compressor.

\begin{figure}[H]
\centering
 \includegraphics[width=0.49\textwidth]{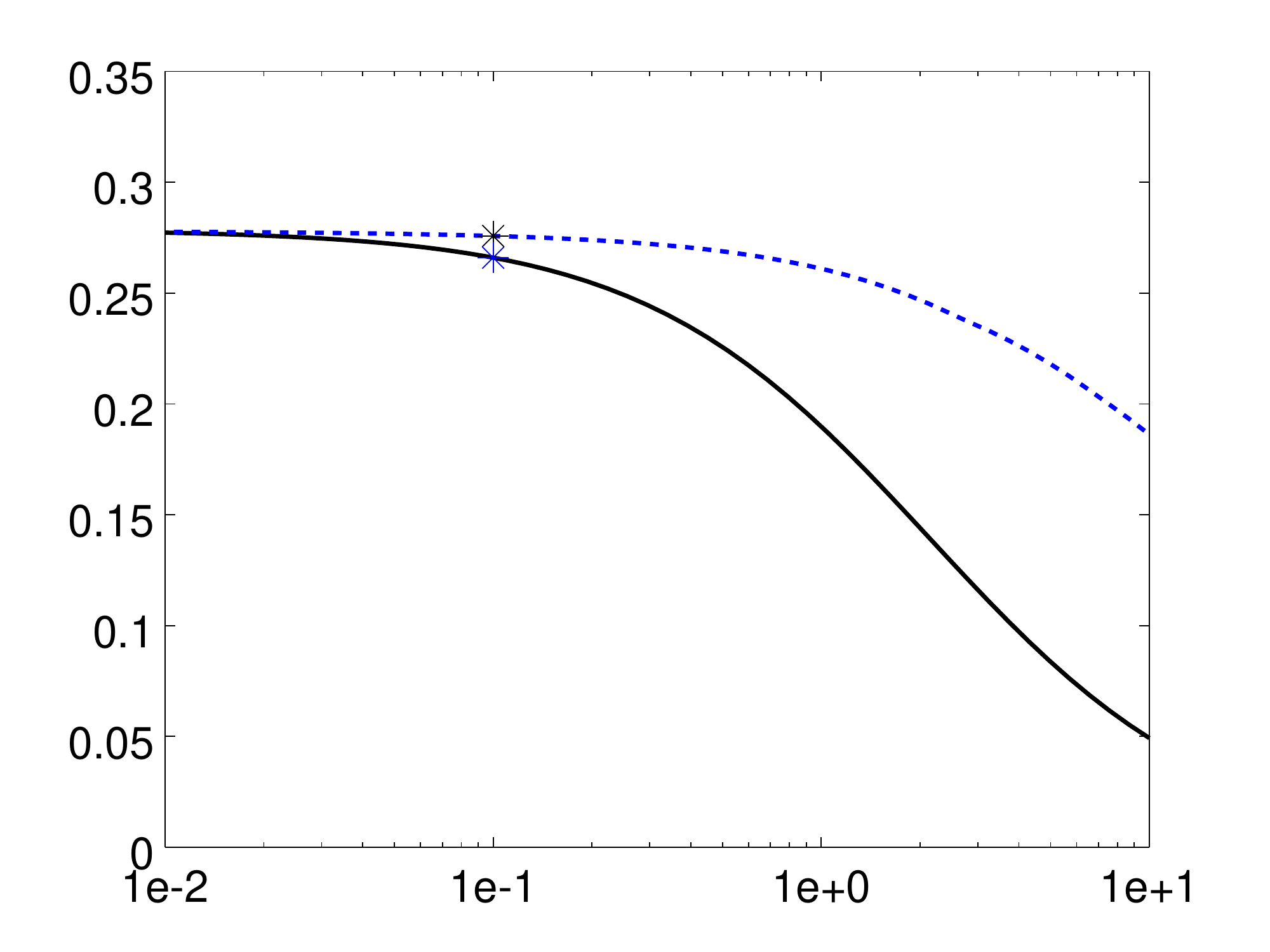}
 \includegraphics[width=0.49\textwidth]{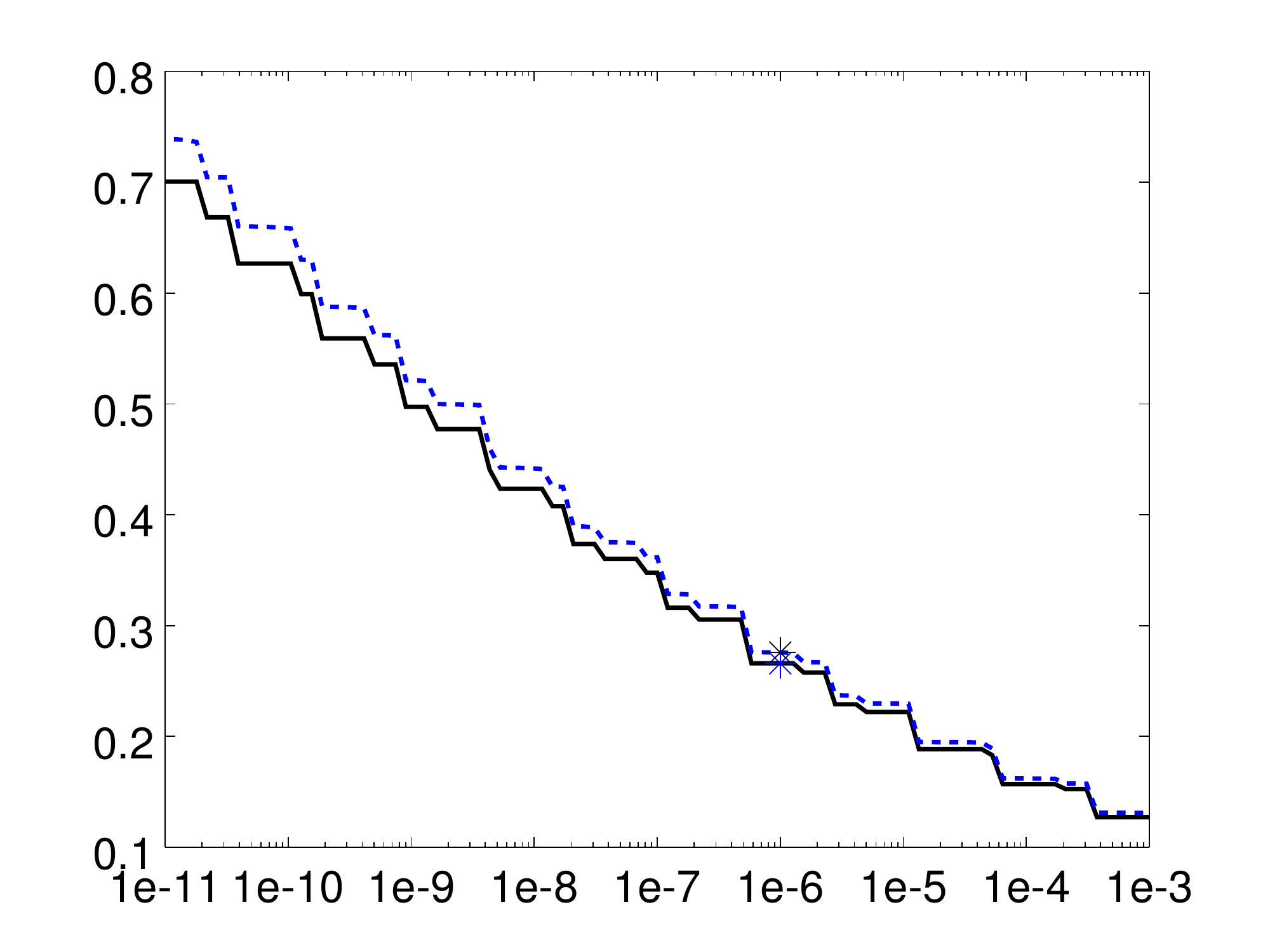}
 \begin{picture}(0,0)(350,0)
  \put(0,110){$E$}
  \put(180,110){$E$}
  \put(80,-5){$t_C(1)$}
  \put(260,-5){TOL}
  \put(90,95){\color{blue} compressed}
  \put(65,50){uncompressed}
  \put(215,100){\color{blue} compressed}
  \put(210,35){uncompressed}
 \end{picture}
 \caption{Theoretically estimated parallel efficiency $E=T_{\rm seq} / (N T_{\rm par})$ for variation of different parameters around the nominal scenario (marked by *). \emph{Left:} varying communication bandwidth in terms of the communication time for uncompressed data. \emph{Right:} varying requested tolerance.}
 \label{fig:parallel-efficiency}
\end{figure}
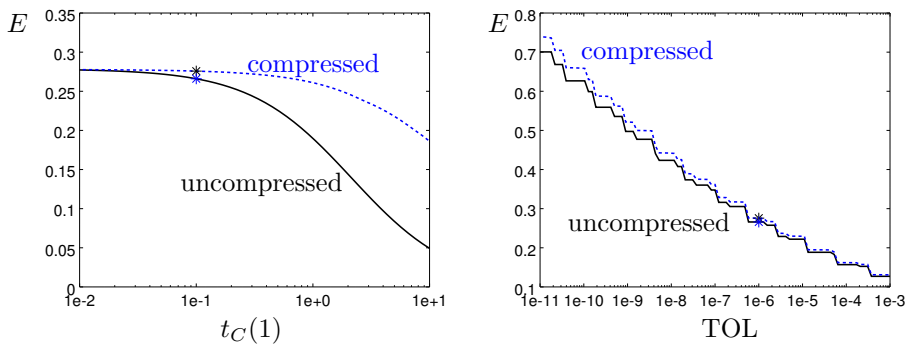

Variation of different parameters in this model around a nominal scenario of contemporary compute clusters as shown in Figure~\ref{fig:parallel-efficiency} suggest that in many current HPC situations, the expected benefit for the run time is small. The predicted run time improvement for the nominal scenario is about 5\%, and does not vary much with, for example, the requested final accuracy (Figure~\ref{fig:parallel-efficiency}, right). The pronounced dependence on smaller bandwidth shown in Figure~\ref{fig:parallel-efficiency}, left, however, makes this approach interesting for a growing imbalance of compute power and bandwidth. Situations where this is already the case is in compute clusters with commodity network hardware and HPC systems where the communication network is nearly saturated due to concurrent communication going on, for example due to the use of spatial domain decomposition.

Indeed, using the cheap transform coding discussed in Section~\ref{subsec:transformcoding} below, an overall run time reduction of 10\% has been observed on contemporary compute nodes connected by gigabit Ethernet~\cite{FischerGoetschelWs2018}.

\subsection{Multilevel Transform Coding on Unstructured Grids for Compressed Communication}\label{subsec:transformcoding}
Due to the larger gap between computing power and bandwidth, transform coding is more attractive for compressing communication in distributed systems than for in-memory compression. While methods based on Cartesian grid structures can be used for some computations, many finite element computations are performed on unstructured grids that do not exhibit the regular tensor structure exploited for designing an orthogonal transform.

An unstructured conforming simplicial grid covers the computational domain $\Omega\subset\R^d$ with non-overlapping simplices $T_i \in \mathcal{T}$, the corners of which meet in the grid vertices $\mathcal{N} = \{x_i\mid i=1,\dots,m\}$; see Figure~9
for a 2D example. The simplest finite element discretization is then with piecewise linear functions, i.e., the solution is sought in the space $V_h = \{ u\in C^0(\Omega) \mid \forall T\in \mathcal{T}: u|_T \in \Pol_1 \}$. The ubiquitous basis for $V_h$ is the nodal basis $(\varphi_i)_{i=1,\dots,m}$ with the Lagrangian interpolation property $\varphi_i (x_j) = \delta_{ij}$, which makes all computations local and leads to sparse matrices.

\begin{figure}[H]
\centering
 \includegraphics[width=\textwidth]{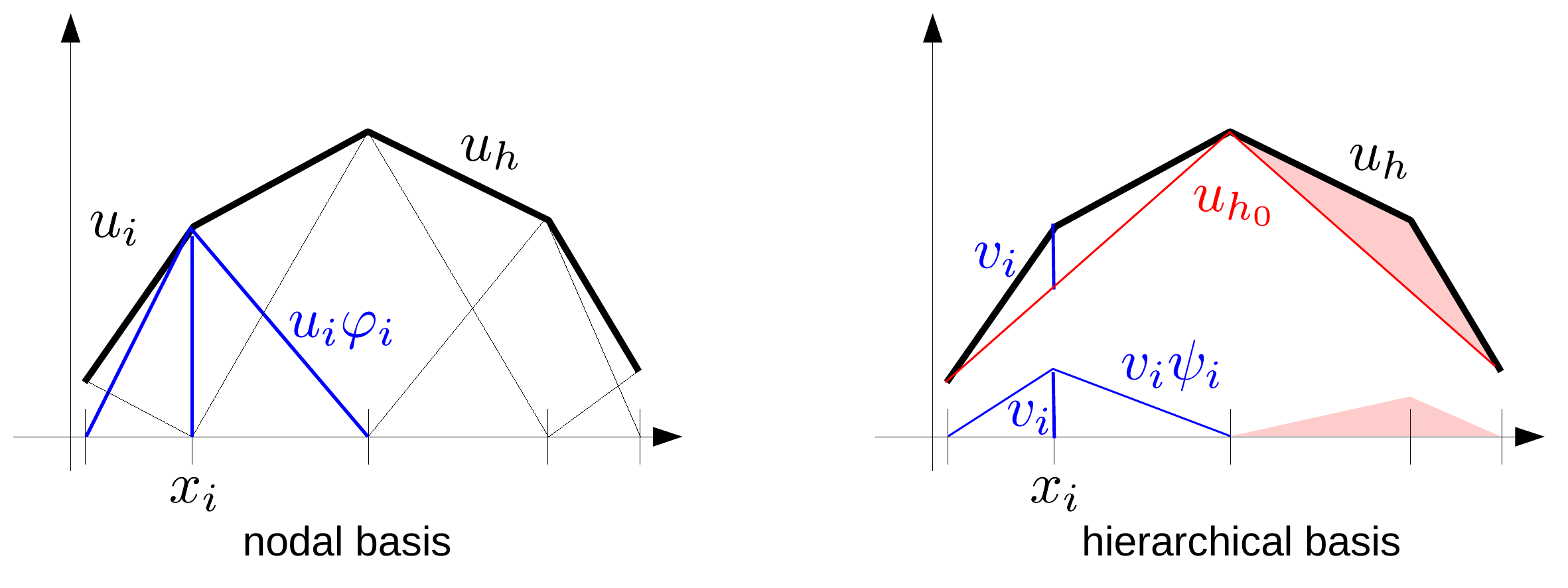}
 \caption{Representation of  1D linear finite element functions $u_h$ in the nodal and hierarchical basis.}
 \label{fig:HB}
\end{figure}

The drawback of the nodal basis is that elliptic systems then lead to ill-conditioned matrices and slow convergence of Krylov methods. Many finite element codes therefore use hierarchies of $\ell+1$ nested grids, resulting from adaptive mesh refinement, for efficient multilevel solvers~\cite{DeuflhardWeiser2012}. The restriction and prolongation operators implemented for those solvers realize a frequency decomposition of the solution
\begin{equation}
 u_h = \sum_{l=0}^{\ell} u_{h_l},
\end{equation}
see Figure~\ref{fig:HB} for a 1D illustration. Using the necessary subset of the nodal basis on grid level $l$ for representing $u_{h_l}$ leads to the hierarchical basis. This hierarchical basis transform allows an efficient in-place computation of optimal complexity and with low overhead, and is readily available in many finite element codes. The transform coefficients can then be quantized uniformly according to the required accuracy and entropy coded~\cite{WsGoetschel2012}, e.g., using a range coder~\cite{Martin79}. Typically, this transform coding scheme (TCUG) takes much less than 5\% of the iterative solution time, see~\cite{WsGoetschel2012,GoetschelVTycowiczPolthierWeiser2015, GoetschelWs2015}.

\begin{figure}[H]
\centering
\scalebox{0.9}{\input{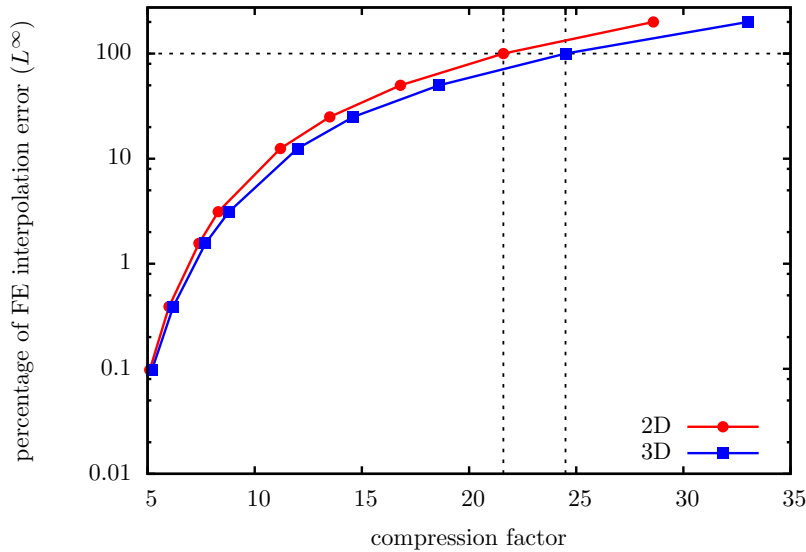}}
\caption[A-priori estimates]{Error vs.~compression factor: A priori estimates for transform coding of finite element functions with hierarchical basis transform, cf.~\cite{WsGoetschel2012}.}
\label{fig:apriorierror}
\end{figure}

A priori error estimates for compression factors and induced distortion can be derived for functions in Lebesgue or Sobolev spaces. The analysis in~\cite{WsGoetschel2012} shows that asymptotically 2.96 bits/value (in 2D, compression factor $21.6$ compared to double precision) are sufficient to achieve a reconstruction error equal to $L^\infty$-interpolation error bounds for functions with sufficient regularity, as is common in elliptic and parabolic equations. In 3D, the compression factor is slightly higher; see Figure~\ref{fig:apriorierror}.

\subsection{Error Metrics} 
An important aspect of compressor design is the norm in which to measure compression errors. While in some applications pointwise error bounds are important and the $L^\infty$-norm is appropriate, other applications have different requirements. For example, if the inexact parallel-in-time method sketched above is applied to parabolic equations, spatially high-frequency error components are quickly damped out. There, the appropriate measure of error is the $H^{-1}$-norm.

Nearly optimal compression factors for given $H^{-1}$-distortion can be achieved in TCUG by replacing the hierarchical basis transform with a wavelet transform, which can efficiently be realized by lifting~\cite{Sweldens1998} on unstructured mesh hierarchies. Level-dependent quantization can be used for near-optimal compression factors matching a prescribed reconstruction error in $H^s$. Rigorous theoretical norm-equivalence results are available for $|s| < 3/2$ with a rather sophisticated construction~\cite{Stevenson2003}. A simpler finite element wavelet construction yields norm equivalences for $-0.114 < s < 3/2$~\cite{CohenEtAl2000},  but in numerical practice it works perfectly well also for $s=-1$. A potential further improvement could be achieved by using rate-distortion theory for allocating quantization levels, as has been done for compression of quality scores in genomic sequencing data~\cite{OchoAsnaniBharadiaChowdhuryWeissmanYona2013}.

\begin{figure}[H]
\centering
\includegraphics[width=0.4\textwidth]{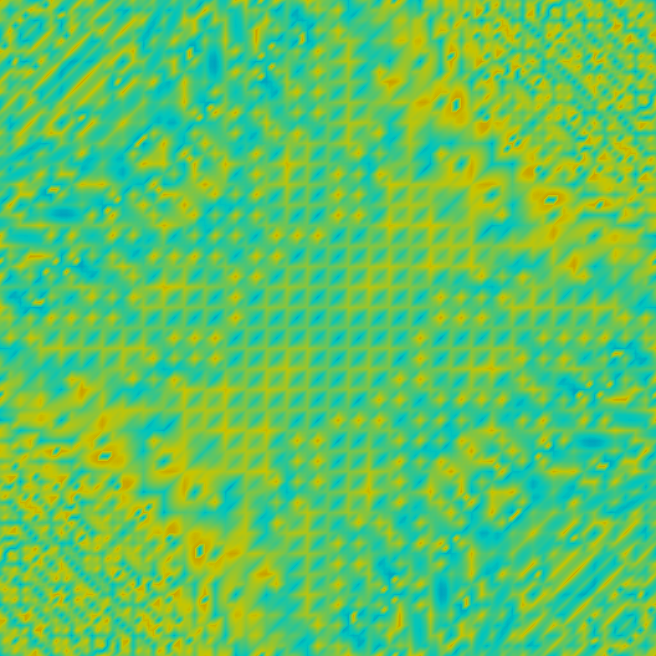}\hspace{2mm}
\includegraphics[width=0.4\textwidth]{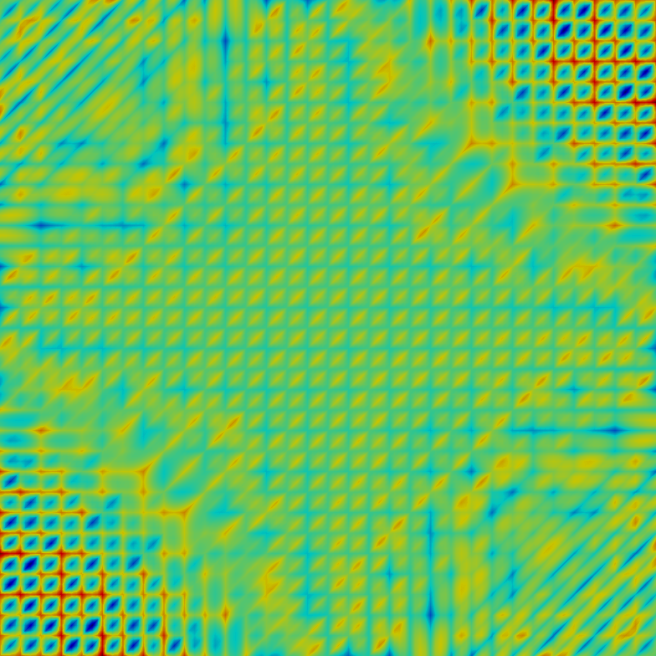}\hspace{2mm}
\includegraphics[height=0.4\textwidth]{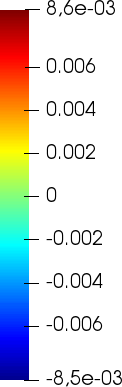}
\caption[Comparison of quantization errors]{Comparison of quantization errors, i.e., error in the reconstructed solution yielding the same $H^{-1}$ error norm. Left: hierarchical basis. Right: wavelets. Using the $H^{-1}$ norm to measure the error allows larger pointwise absolute reconstruction errors compared to the $L^\infty$ error metric, thus higher compression factors.}
\label{fig:quanterr}
\end{figure}

Figure~\ref{fig:quanterr}  shows the quantization errors for the 2D test function $f(x) = \sin(12(x_0-0.5)(x_1-0.5))$ on a uniform mesh of $16{,}641$ nodes, with a grid hierarchy of seven levels. Using a wavelet transform almost doubles the compression factor here, while keeping the same $H^{-1}$ error bound as the hierarchical basis transform~\cite{Goe2015}.

A closely related aspect is the order of quantization and transform. In the considerations above, a transform-then-quantize
 approach has been assumed. An alternative is the quantize-then-transform~sequence, which then employs an integer transform. It allows guaranteeing strict pointwise reconstruction error bounds directly, and is therefore closely linked to $L^\infty$ error concepts. In contrast, quantization errors of several hierarchical basis or wavelet coefficients affect a single point, i.e., a single nodal basis coefficient. The drawback of quantize-then-transform is that the quantization step shifts energy from low-frequency levels to high-frequency levels, leading to less efficient decorrelation if error bounds in Sobolev spaces are important. A brief discussion and numerical comparison of the two approaches for some test functions can be found in~\cite{WsGoetschel2012}.

\section{Mass Storage} \label{sec:storage}

The third level in the compression hierarchy is mass storage. Often, single hard disks or complete storage systems are again slower than inter-node communication links in distributed memory systems. The significantly higher flops/byte ratio makes more sophisticated and more effective compression schemes attractive, and in particular allows employing a posteriori error estimators for a better control of the compression error tolerance. These schemes are necessarily application-specific, since they need to predict error transport into the final result, and to anticipate the intended use of reconstructions as well as the required accuracy. Examples are adjoint gradient computation in PDE-constrained optimization problems (Section~\ref{sec:adjoint}) and checkpoint/restart for fault tolerance (Section~\ref{sec:checkpoint}).

In the extreme case, the size of the data to store is the limiting factor, and the computational effort for compression does not play a significant role. This is typically the case in solution archiving (Section~\ref{sec:archiving}).

\subsection{Adjoint Solutions}\label{sec:adjoint}
% \cite{GoetschelChamakuriKunischWeiser2014,GoetschelWs2015,GoetschelTycowiczPolthierWeiser2015,WsGoetschel2012}

Adjoint, or dual, equations are important in PDE-constrained optimization problems, e.g., optimal control in electrophysiology~\cite{GoetschelChamakuriKunischWeiser2014} or inverse problems~\cite{BoeHanPueFic2016,LinCheLee2016}, and goal-oriented error estimation~\cite{OdePru2001}. Consider the abstract variational problem
\begin{equation}\label{eq:abstract_var}
\text{find}\ x\in X\ \text{such that}\ c(x;\varphi) = 0 \ \forall \varphi \in \Phi,
\end{equation}
for a differentiable semilinear form $c : X \times \Phi \rightarrow Z$ with suitable function spaces $X,\Phi, Z$, and a quantity of interest given as a functional $J : X \rightarrow \R$. During the numerical solution, Equation~(\ref{eq:abstract_var}) is typically only fulfilled up to a nonzero residual $r$, i.e., $c(x;\varphi) = r$. Naturally the question arises, how does the residual $r$ influence the quantity of interest $J(x)$. For instationary PDEs, answering this question leads to solving an adjoint equation backwards-in-time. As the adjoint operator and/or right-hand sides depend on the solution $x$, storage of the complete trajectory is needed, thus requiring techniques to reduce the enormous storage demand for large-scale, real-life applications. We note in passing, that, obviously, compression is not only useful for storage on disk, but can also be used in-memory, thus allowing  more data to be kept in RAM and potentially avoiding disk access.

Lossy compression for computing adjoints can be done using transform coding as discussed in Section~\ref{subsec:transformcoding}. In addition to the spatial smoothness, correlations in time can be exploited for compression. Since the stored values are only accessed backwards in time, following the adjoint equation integration direction, it is sufficient to store the state at the final time and its differences between successive time~steps. This predictive coding with a constant state model, also known as delta encoding, can be efficiently implemented, requiring only to keep one additional time step in memory. Linear or higher order models can be used for prediction in time as well, but already the most simple delta encoding can significantly increase---in some cases double---the compression factor at very small computational cost. For more details and numerical examples we refer to~\cite{WsGoetschel2012,GoetschelVTycowiczPolthierWeiser2015}.

Before presenting examples using lossy compression for PDE-constrained optimization and goal-oriented error estimation, let us briefly mention so-called checkpointing methods for data reduction in adjoint computations, first introduced by Volin and Ostrovskii~\cite{VolOst1985}, and Griewank~\cite{Gri1992}. Instead of keeping track of the whole forward trajectory, only the solution at some intermediate timestep is stored. During the integration of the adjoint equation, the required states are re-computed, starting from the snapshots, see, e.g.,~\cite{GriWal2008} for details. This increases the computational cost for typical settings (compression factors around 20) by two to four additional solves of the primal PDE.  Moreover, due to multiple read- and write-accesses of checkpoints during the re-computations for the adjoint equation, the reduction in memory \emph{bandwidth} requirements is significantly smaller. We refer the reader to~\cite{GoetschelVTycowiczPolthierWeiser2015} for a more detailed discussion and additional references.

\subsubsection{PDE-Constrained Optimization} 
For PDE-constrained optimization, typically $X = Y \times U, x = (y,u)$ in the abstract problem~\eqref{eq:abstract_var}, where the influence of the control $u$ on the state $y$ is given by the PDE. Here, $J$ is the objective to be minimized, e.g., penalizing the deviation of $y$ from some desired state. Especially in time-dependent problems, often the reduced form is considered: There, the PDE~\eqref{eq:abstract_var} is used to compute for a given control $u$ the associated (locally) unique solution $y=y(u)$. With only the control remaining as the optimization variable, the reduced problem reads $\min_u j(u)$, with $j(u):=J(y(u),u)$. Computation of the reduced gradient then leads to the adjoint equation for $p\in Z^\star$
\begin{equation}\label{eq:adj_grad}
c_y^\star(p;(y,u),\varphi) = -J_y((y,u),\varphi),
\end{equation}
where $\star$ denotes the dual operator/dual function spaces, and $c_y, J_y$ are the derivatives of $c(y,u;\varphi), J(y,u)$ with respect to the $y$-component.

Exemplarily, we consider optimal control of the monodomain equations on a simple 2D unit square domain $\Omega$.
This system describes the electrical activity of the heart (see, e.g.,~\cite{ColDdErdLanPav2006}), and consists of a pa\-ra\-bo\-lic PDE for the transmembrane voltage $v$, coupled to pointwise ODEs for the
gating variable $w$ describing the state of ion channels:
\begin{equation}\label{eq:monodomain}
\begin{alignedat}{3} 
 v_t &= \mathop\mathrm{div}( \sigma \nabla v) - I_\text{ion}(v,w)  + I_\text{e} &\quad& \textrm{in}~\Omega \times (0,T)\\
 w_t &= G(v,w) && \textrm{in}~\Omega \times (0,T).
\end{alignedat}
\end{equation}
Homogeneous Neumann boundary conditions are prescribed. The functions $I_\text{ion}(v,w)$ and $G(v,w)$ are specified by choosing a membrane model.
%\begin{align*}
%I_\text{ion}(v,w) &= g v \bigl(1-\frac{v}{v_{th}}\bigr)\bigl(1-\frac{v}{v_{p}}\bigr) + \eta_1 v w\\
%G(v,w) &= \eta_2 \bigl( \frac{v}{v_p} - \eta_3 w\bigr)
%\end{align*}
%and homogeneous Neumann boundary conditions. % and initial values.
%In this 2D model $\sigma : \R^2 \rightarrow \R^{2\times2}$ and $g, \eta_i, v_p, v_{th} \in \R_+$ are given parameters (see, e.g.,~\cite{ColDdErdLanPav2006}). %Table~\ref{tab:electro-parameters}). 
%\begin{table}[ht]
%\begin{tabular}{cccccccc}
% \toprule
% $\sigma_{il}$ & $\sigma_{it}$ & $g$ & $v_{th}$ & $v_p$ & $\eta_1$ & $\eta_2$ & $\eta_3$ \\
% $[\Omega^{-1} \mathrm{cm}^{-1}]$ &  $[\Omega^{-1} \mathrm{cm}^{-1}]$ &$[\mathrm{mS}/\mathrm{cm}^2]$ & $[\mathrm{mV}]$ & $[\mathrm{mV}]$ & $[\mathrm{mS}/\mathrm{cm}^2]$ & &
% \\ \midrule
% $3\cdot 10^{-3}$ & $3.1525\cdot 10^{-4}$ & $1.5$ & $ 13$ &$100$ & $4.4$ & $0.012$ & $1$
% \\ \bottomrule
%\end{tabular}
%\caption[Electrophysiological parameters]{Electrophysiological parameters (adapted from~\cite{ColDdErdLanPav2006}).}\label{tab:electro-parameters}
%\end{table}
For the optimal control problem an initial excitation in some subdomain $ \Omega_\text{exi}$ is prescribed. %initial values
%\begin{align*}
%   v(x,0) &= \begin{cases}101.0 \quad &\text{in} \quad \Omega_\text{exi} \\
%  0 &\text{otherwise} \end{cases} \nonumber \\
%  w(x,0) &= 0 \quad \text{in}~\Omega.
%\end{align*}
%are prescribed.
%Here, $\Omega_\text{exi}$ is a circle with midpoint $(0.5, 0.5)$ and radius $0.04$. 
The external current stimulus is $I_\text{e}(x,t) = \chi_{\Omega_c}(x) u(t)$, where the control $u$ is spatially constant on a control domain $\Omega_{c}$. % =[0.37,0.4]\times [0.45,0.55]\cup [0.6,0.63]\times [0.45,0.55]$.
Defining some observation domain $\Omega_\text{obs}$, the objective is given by
\begin{equation}\label{eq:mono_obj}
  J(y,u) = \half\norm{v}^2_{L^2(\Omega_\text{obs} \times (0,T))} + \frac{\alpha}{2} \norm{u}^2_{L^2(0,T)},
\end{equation}
i.e.,~we aim at damping out the excitation wave. %We use
%\begin{equation}
%\Omega_\text{obs} = \Omega \setminus \bigl( [0.35,0.42]\times [0.43,0.57] \cup [0.58,0.65]\times [0.43,0.57] \bigr), 
%\end{equation}
% $T=4$, and $\alpha = 3\cdot 10^{-6}$;
For details, see~\cite{NagaiahKunischPlank2011}.
Solution of the optimization problem with inexact Newton-CG methods and lossy compression is investigated in~\cite{GoetschelChamakuriKunischWeiser2014,GoetschelWs2015}; here we use the quasi-New\-ton meth\-od due to Broyden, Fletcher, Goldfarb, and Shanno (BFGS, see, e.g.,~\cite{DennisMore1977,BorziSchulz2012}). For time discretization, we use a linearly implicit Euler method with fixed timestep size $dt=0.04$. Using linear finite elements, spatial adaptivity is performed individually for state and adjoint using a hierarchical error estimator~\cite{DdLeiYse1989}, with a restriction to at most $25,000$ vertices in space. The adaptively refined grids were stored using the methods from~\cite{TycowiczKaelbererPolthier2011}, which reduced the storage space for the mesh to less than 1 bit/vertex (see~\cite{GoetschelVTycowiczPolthierWeiser2015}).

Lossy compression of state values, i.e., the finite element solutions $v$ and $w$, at all time steps, affects the accuracy of the reduced gradient computed by adjoint methods, and results in inexact quasi-Newton updates. Error analysis~\cite{Goe2015} shows that BFGS with inexact gradients converges linearly, if the gradient error $e_g$ in each step fulfills
 \begin{equation}\label{eq:QNerrorTol}
  \norm{e_g} \leq \frac{\varepsilon}{\kappa(B)^{1/2}} \norm{\tilde{g}}
 \end{equation}
 for $\varepsilon < 1/2$. Here, $\kappa(B)$ is the condition number of the approximate Hessian $B$, and $\tilde{g}$ denotes the inexactly computed gradient. The error bound~\eqref{eq:QNerrorTol} allows computing adaptive compression error tolerances from pre-computed worst-case gradient error estimates analogously to~\cite{GoetschelWs2015}, see~\cite{Goe2015} for~details.
 
\begin{figure}[H]
\centering
\scalebox{0.85}{\input{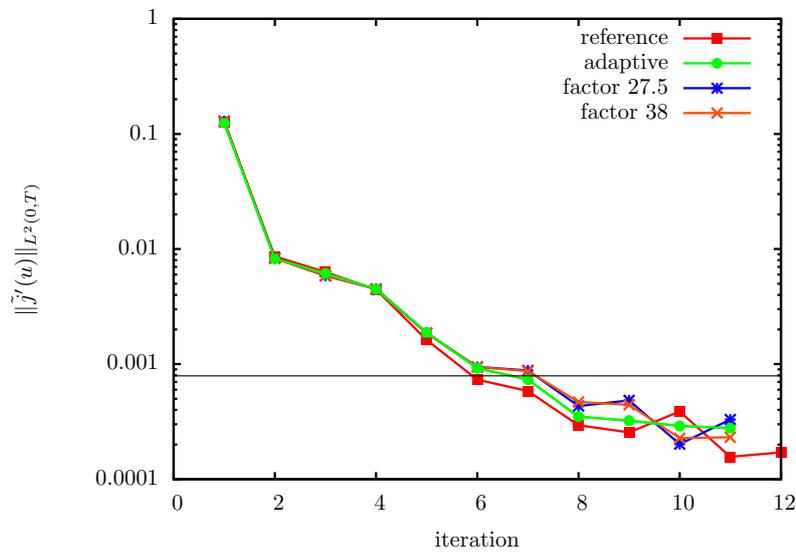}}
\caption[Optimization progress for BFGS (monodomain example)]{Optimization progress of Broyden, Fletcher, Goldfarb, and Shanno (BFGS) for the monodomain example~\eqref{eq:monodomain}, \eqref{eq:mono_obj}, using different quantization tolerances for the state trajectory. No delta-encoding between timesteps was used. The horizontal line shows the approximate discretization error of the reduced gradient. See also~\cite{GoetschelVTycowiczPolthierWeiser2015}.}
\label{fig:BFGSprogress}
\end{figure}

\begin{figure}[H]
% \begin{minipage}{\textwidth}
\centering
% \makebox[\textwidth]{%
\scalebox{0.85}{\input{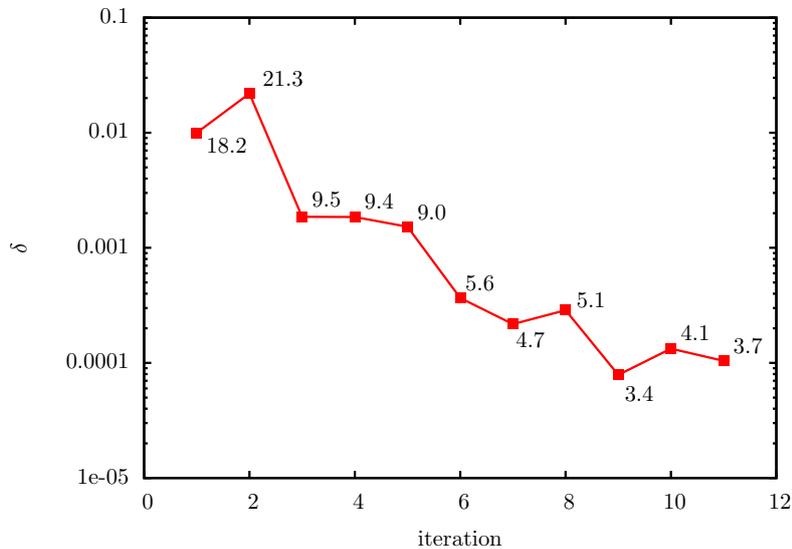}}%
% }
\caption[Adaptive quantization tolerances for BFGS (monodomain example)]{Adaptively chosen quantization tolerances $\delta$ and cor\-re\-spond\-ing compression factors for the monodomain example~\eqref{eq:monodomain},~\eqref{eq:mono_obj} using BFGS. The relatively small adaptive compression factors are due to using worst case error estimates in~\eqref{eq:QNerrorTol}, as well as a fixed maximum mesh size.}
\label{fig:BFGSfactors}
% \end{minipage}
\end{figure}

 Figure~\ref{fig:BFGSprogress} shows the progress of the optimization method. For trajectory compression, different fixed as well as the adaptively chosen quantization tolerances were used. We estimate the spatial discretization error in the reduced gradient by using a solution on a finer mesh as a reference.
Up to discretization error accuracy, lossy compression has no significant impact on the optimization progress.
The adaptively chosen quantization tolerances for the state values are shown in Figure~\ref{fig:BFGSfactors}.  The resulting compression factors when using TCUG as discussed in Section~\ref{subsec:transformcoding} together with delta encoding in time shown as well. In the first iteration, a user-prescribed tolerance was used. %Due to the relatively large reduction of the gradient norm in the first iteration, the gradient norm estimate by equation~\eqref{c5:eq:theta} used for determination of $\delta^y$ for iteration 3 is too strict, such that the chosen tolerance is too small.
The small compression factors are on one hand due to the compression on adaptively refined grids, as discussed in Section~\ref{sec:adaptive}, and on the other hand due to overestimation of the error in the worst case error estimates. The latter is apparent from the comparison with prescribed fixed quantization tolerances (see Figure~\ref{fig:BFGSprogress}). To increase the performance of the adaptive method, tighter, cheaply computable error estimates are required.

% The estimated condition number of the reduced Hessian varies between $200$ and $230$. We note that compared to prescribed fixed quantization error tolerances, the adaptively chosen tolerances are too restrictive (see Fig.~\ref{fig:BFGSprogress}). This is , and the fixed tolerance for the discretization.

\subsubsection{Goal-Oriented Error Estimation}
For goal-oriented adaptivity, we consider solving the PDE~\eqref{eq:abstract_var} by  a Galerkin approximation,
\begin{equation}\label{eq:abstract_var_Gal}
\text{find}\ x_h\in X_h\ \text{such that}\ c(x_h;\varphi_h) = 0 \quad \forall \varphi_h \in \Phi_h,
\end{equation}
with suitable finite dimensional subspaces $X_h\subset X, \Phi_h \subset \Phi, Z_h \subset Z$. Here the functional $J$ measures some quantity of interest, e.g., the solution's value at a certain point, or in case of optimal control problems the objective, with the aim that
\begin{equation}
|J(x_h) - J(x)| \leq \epsilon.
\end{equation}
The dual weighted residual (DWR) method~\cite{BecRan1996,BecKapRan2000} now seeks to refine the mesh used to discretize the PDE by weighting (local) residuals with information about their global influence on the goal functional $J$~\cite{Ran2006}. These weights are computed by the dual problem
\begin{equation}\label{eq:dwr_dual}
\text{find}\ p \in Z^\star\ \text{such that}\ c_x^\star(p;x_h,\varphi) = J_x(x_h,\varphi) \quad \forall \varphi \in \Phi^\star,
\end{equation}
which depends on the approximate primal solution $x_h$, which therefore needs to be stored.

\begin{figure}[H]
\centering
	\includegraphics[width=0.235\textwidth]{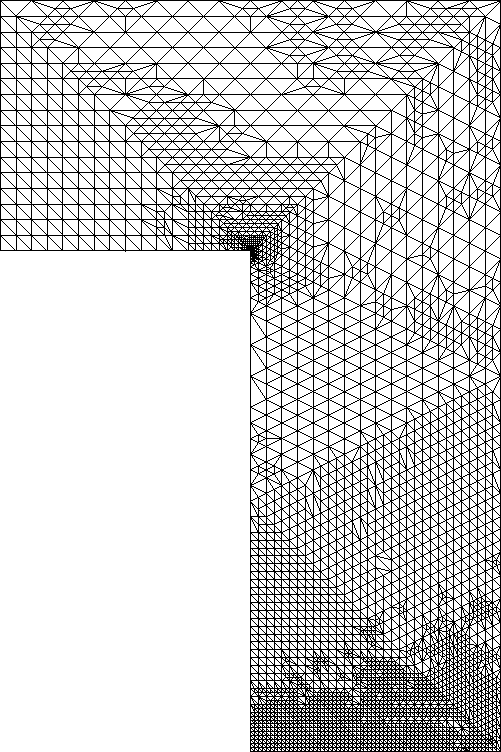}% RESOLVED: "line" break in figure
	\includegraphics[width=0.235\textwidth]{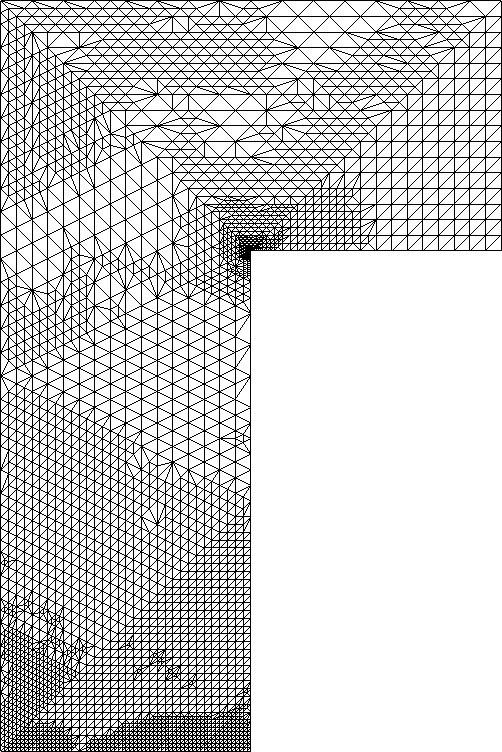}%
	\includegraphics[width=0.525\textwidth]{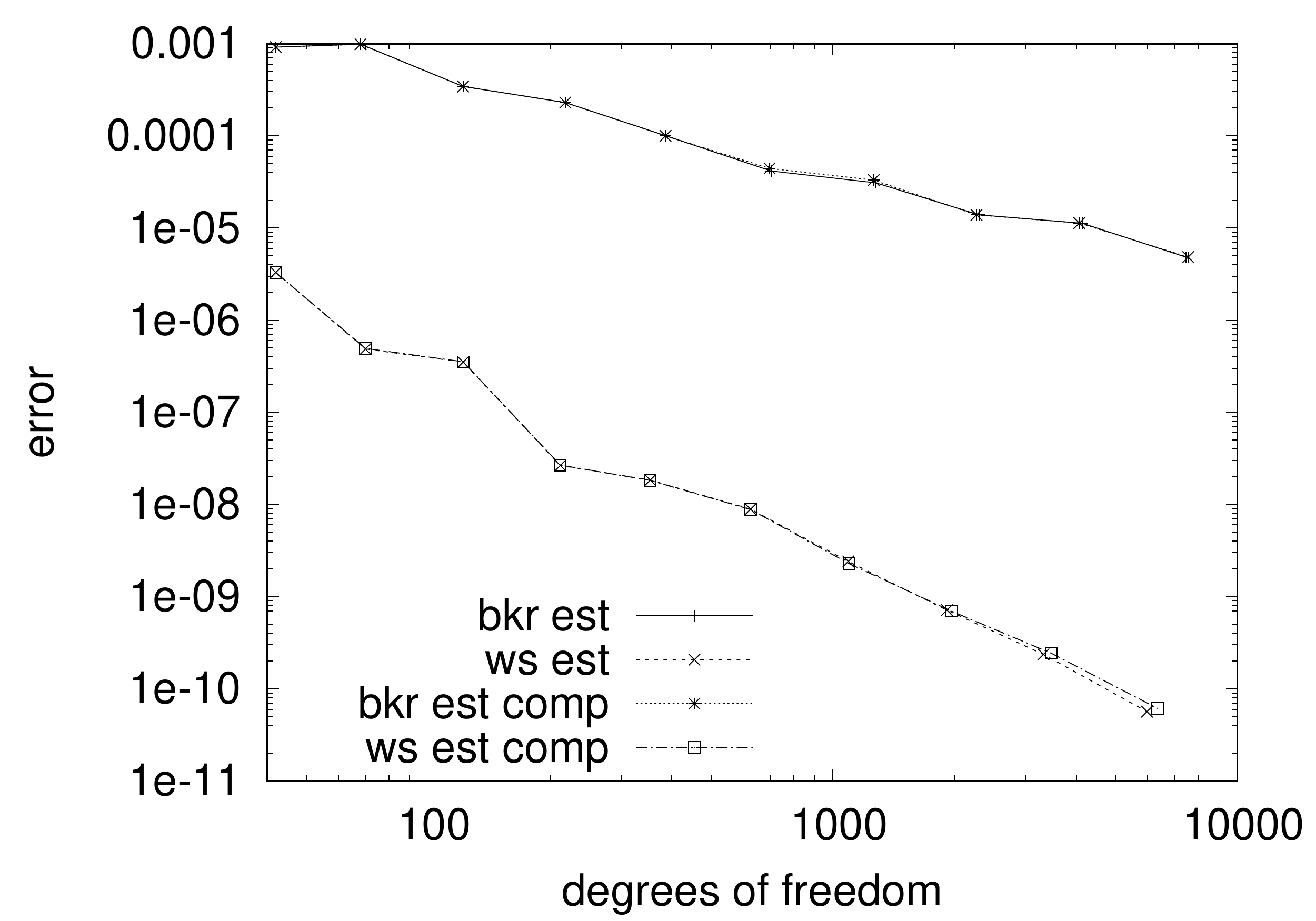}
	\caption{Goal-oriented error estimation: generated meshes with (left) and without (middle) compression (with compression factor up to 32) for Example 3b in~\cite{Ws2013}. Meshes were generated using weights according to Weiser~\cite{Ws2013}. Estimated errors (right) are shown for weights due~\cite{Ws2013} (ws) and to Becker et al.~\cite{BecKapRan2000} (bkr), both with and without compression. Differences between estimated errors with and without compression are barely visible, and negligible compared to the differences between the two error concepts.}
	\label{fig:errorest}
\end{figure}

The effect of compression on error estimation is illustrated in Figure~\ref{fig:errorest}. For simplicity we use a linear-quadratic elliptic optimal control problem here (Example 3b in~\cite{Ws2013}), with the objective as goal functional. Extension to time-dependent problems using the method of time layers (Rothe method) for time discretization is straightforward. Meshes were generated using weights, according to Weiser~\cite{Ws2013}, for estimating the error in the reduced functional $J(y(u_h),u_h)-J(\bar y, \bar u)$, as well as due to \mbox{Becker et al.~\cite{BecKapRan2000}} for the all-at-once error $J(y_h, u_h)-J(\bar y, \bar u)$. The compression tolerance for TCUG was chosen such that the error estimation is barely influenced, resulting in only slightly different meshes. For the final refinement step, i.e., on the finest mesh, compression resulted in data reduction by a factor of 32.
Thus, at this factor the impact of lossy compression is barely noticeable in the estimated error \mbox{as well as} in the resulting mesh, and is much smaller than, e.g., the influence of the chosen error concept for mesh refinement. Numerical experiments using various compression techniques can be found in~\cite{CyrShadWil2015}. Instead of the ad-hoc choice of compression tolerances, a thorough analysis of the influence of compression error on the error estimators is desirable; this is, however, left for future work. 

\subsubsection{Adaptive Grid Refinement} \label{sec:adaptive}

PDE solutions often exhibit spatially local features, e.g., corner singularities or moving fronts as in the monodomain Equation~(\ref{eq:monodomain}), which need to be resolved with small mesh width. Uniform grids with small mesh width lead to huge numbers of degrees of freedom, and therefore waste computing effort, bandwidth, and storage capacity in regions where the solution is smooth. Adaptive mesh refinement, based on error estimators and local mesh refinement, has been established as an~efficient means to reduce problem size and solution time, cf.~\cite{DeuflhardWeiser2012}, and acts at the same time as a decimation/interpolation-based method for lossy data compression. The construction of adaptive meshes has been explicitly exploited by Demaret et al.~\cite{DemaretDynFloaterIske2005} and Solin and Andrs~\cite{SolinAndrs2009} for scientific data~compression.

\begin{figure}[H]
\centering
 \includegraphics[width=0.55\textwidth]{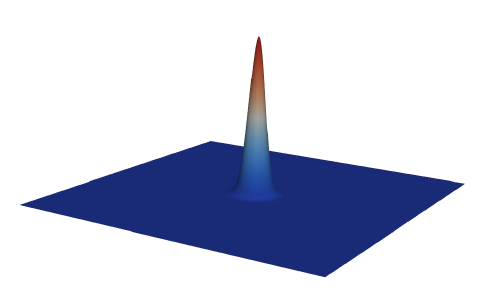}
 \includegraphics[width=0.35\textwidth]{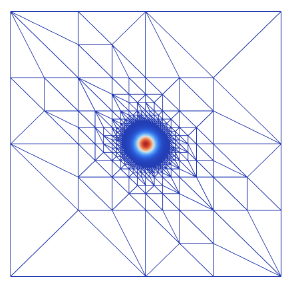}
 \caption{\emph{Left:} highly local peak function. \emph{Right:} adaptively refined mesh with 4237 vertices for minimal finite element interpolation error. A uniform grid with the same local resolution has 263,169~vertices.}
 \label{fig:adaptive}
\end{figure}

Interestingly, even though decimation and hierarchical transform coding both compete for the same spatial correlation of data, i.e., smoothness of functions to compress, their combination in PDE solvers can achieve better compression factors for a given distortion than each of the approaches~alone. A simple example is shown in Figure~\ref{fig:adaptive} with compression factors for a fixed error tolerance given in Table~\ref{tab:adaptivity}. Using both, adaptive mesh refinement and transform coding, is below the product of individual compression factors, which indicates that there is in fact some overlap and competition for the same correlation budget. Nevertheless, it shows that even on adaptively refined grids there is a significant potential for data compression. The compression factor of adaptive mesh refinement as given in Table~\ref{tab:adaptivity} is, however, somewhat too optimistic, since it only counts the number of coefficients to be stored. For reconstruction the mesh has to be stored as well. Fortunately, knowledge about the mesh refinement algorithm can be used for extremely efficient compression of the mesh structure~\cite{GoetschelVTycowiczPolthierWeiser2015}.

\begin{table}[H]
\centering
 \caption{Compression factors for adaptive mesh refinement and transform coding of the peak function shown in Figure~\ref{fig:adaptive}. Adaptive and uniform mesh refinement yield the same interpolation error; the same error tolerance for lossy compression was used in the two cases.}
 \label{tab:adaptivity}
 \begin{tabular}{ccc}
 \toprule
           & \textbf{Double} & \textbf{Transform Coding} \\
  \midrule
  uniform  &    1   &  54 \\
  adaptive &   62   & 744 \\
  \bottomrule
 \end{tabular}

\end{table}

Another reason why adaptive mesh refinement and transform coding can be combined effectively for higher compression factors, is that the accuracy requirements for the PDE solution, and the solution storage can be very different, depending on the error propagation.
Thus, the decimation by mesh adaptivity may need to retain information that the transform coder can safely neglect, allowing the latter to achieve additional compression on top of the former. For example, in adjoint gradient computation, the mesh resolution affects the accuracy of all future time steps and via them the reduced gradient, whereas the solution storage for backwards integration of the adjoint equation affects only one time step.

\subsection{Checkpoint/Restart}\label{sec:checkpoint}
In exa-scale HPC systems, node failure will be a common event. Checkpoint/re\-start is thus mandatory, but snapshotting for fault tolerance is increasingly expensive due to checkpoint sizes. Application-based checkpointing aims at reducing the overhead by optimizing snapshot times, i.e., when to write a checkpoint, and what to write, i.e., store only information that cannot be re-computed in a reasonable amount of time. Moreover, information should only be stored with required accuracy, which might be significantly smaller than double precision values.

Lossy compression for fault-tolerant iterative methods to solve large-scale linear systems is discussed by Tao et al.~\cite{TaoDiLiaCheCap2018}. They derive a model for the computational overhead of checkpointing both with and without lossy compression, and analyze the impact of lossy checkpointing. Numerical experiments demonstrate that their lossy checkpointing method can significantly reduce the fault tolerance overhead for the Jacobi, GMRES, and CG methods.

Calhoun et al.~\cite{CalCapOlsSniGro2019} investigate using lossy compression to reduce checkpoint sizes for time stepping codes for PDE simulations. For choosing the compression tolerance they aim at an error less than the simulation's discretization error, which is estimated a priori using information about the mesh width of the space discretization and order of the numerical methods. Compression is performed using SZ~\cite{DiCap2016,TaoDiCheCap2017}. Numerical experiments for two model problems (1D-heat and 1D-advection equations) and two HPC applications (2D Euler equations with PlacComCM, a multiphysics plasma combustion code, and 3D Navier-Stokes flow with the code Nek5000) demonstrate that restart from lossy compressed checkpoints does not significantly impact the simulation, but reduces checkpoint time.

\begin{figure}[H]
\centering
\includegraphics[width=0.45\textwidth]{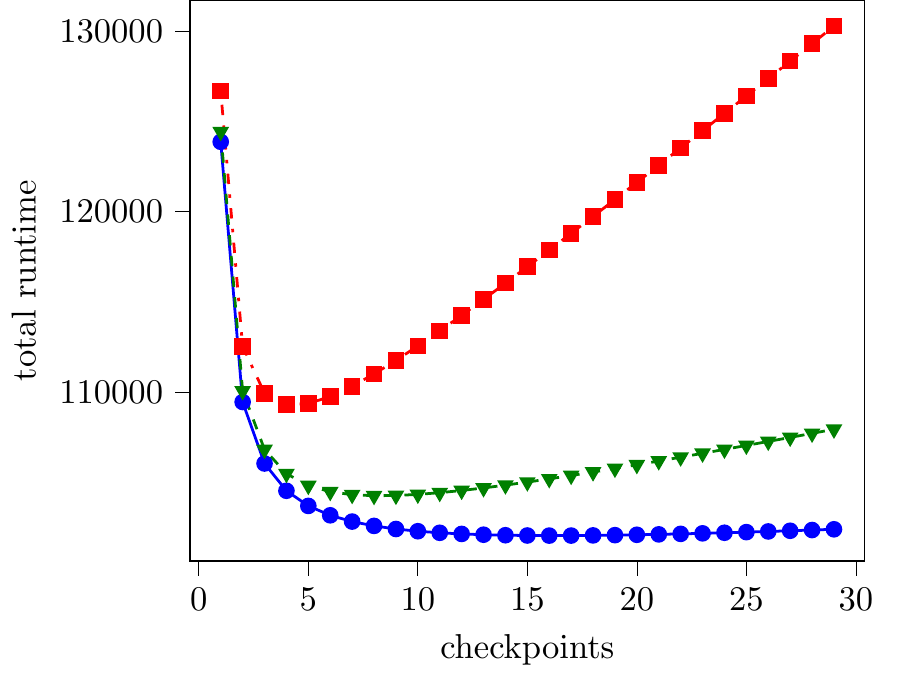}
%\scalebox{0.225}{\input{T_vs_n.tex}}
\includegraphics[width=0.45\textwidth]{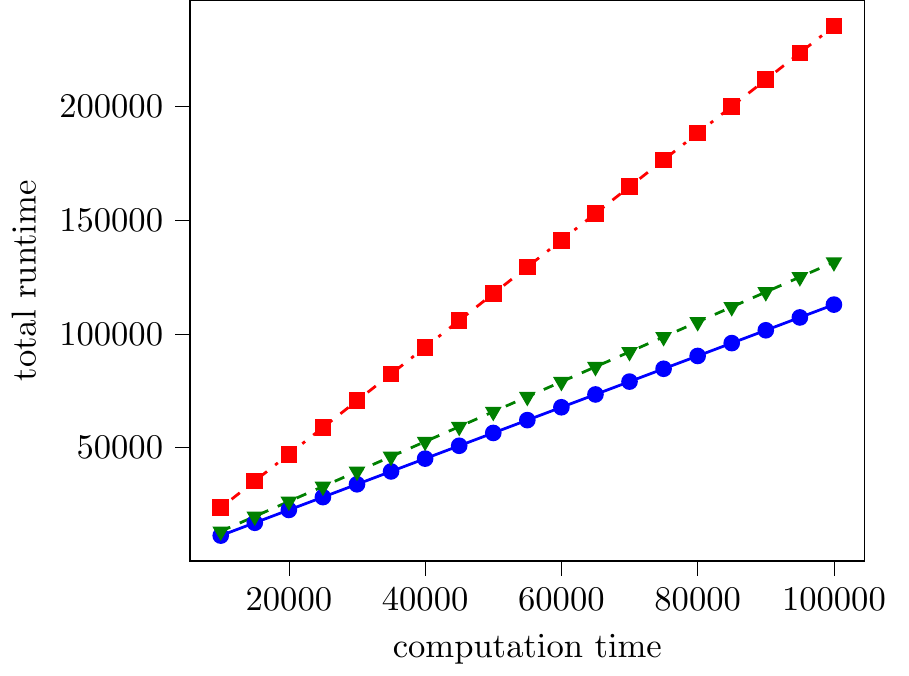}\\
%\scalebox{0.225}{\input{T_vs_Tc.pgf}}\\
\includegraphics[scale=1]{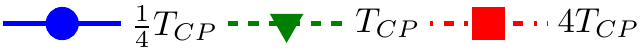}
\caption{Influence of $T_\text{CP}$, comparing the nominal scenario (green) to scenarios with $T_\text{CP}$ reduced (blue) and increased (red) by a factor 4. \emph{Left:} overall runtime vs. number of checkpoints for $N=4$ and $T_C = 100000\mathrm{s}$. \emph{Right:} overall runtime vs. actual computation time for $N=100$. In both cases $p_\text{RS}=7.74\cdot10^{-7}$, $T_\text{CP}=245.583\mathrm{s}$ and $T_R=545.583\mathrm{s}$ were used. While the runtimes were measured for the parallel-in-time solution of a 3D heat equation on the HLRN-III Cray XC30/XC40 supercomputer (\url{www.hlrn.de}), the probability of failure was determined from HLRN-III logfiles.}
\label{fig:checkpointing}
\end{figure}

Application-specific fault-tolerance including computing optimum checkpoint intervals~\cite{Young1974,Daly2006} or multilevel checkpointing techniques~\cite{DiRobVivCap2017} has been a research topic for many years. In order to illustrate the influence of lossy compression, we derive a simple model similar to~\cite{Daly2006}, relating probability of failure and checkpoint times to the overall runtime of the application.  We assume equidistant checkpoints and aim at determining the optimum number of checkpoints $n$. For this we consider a parallel-in-time simulation application (see Section~\ref{sec:PinT}) and use the notation summarized in Table~\ref{tab:checkpointing-notation}.
The overall runtime of the simulation consists of the actual computation time $T_C$, the time it takes to write a checkpoint $T_\text{CP}$, and the restart time $T_\text{RS}$ for $N_\textbf{RS}$ failures/restarts,
 $T = T_C + n T_\text{CP} + T_\text{RS} N_\text{RS}$. Note that here $T_C$ depends implicitly on the number of cores used. The time for restart $T_\text{RS}$ consists of the average required re-computation from the last written checkpoint to the time of failure, here for simplicity assumed as $\frac{1}{2} \frac{T}{n}$ (see also~\cite{Young1974}) and time to recover data structures $T_\text{R}$. 
 For $N$ compute cores, and probability of failure per unit time and core $p_\text{RS}$, we get the estimated number of restarts $N_\text{RS} = p_\text{RS} T N$. Bringing everything together, the overall runtime amounts to
\begin{equation}
T(n) = T_C + n T_\text{CP} + T_\text{RS} N_\text{RS} = \frac{n\bigl(b-\sqrt{b^2-\frac{2}{n}p_\text{RS} N (T_\text{C}+nT_\text{CP})}\bigr)}{p_\text{RS}N},
\end{equation}
where for brevity we use the unit-less quantity $b = 1-T_\text{R}p_\text{RS}N$. As $T_\text{RS}$ and $N_\text{RS}$ depend on the total time $T$, this model includes failures during restart as well as multiple failures during the computation of one segment between checkpoints.
Given parameters of the HPC system and the application, an optimal number of checkpoints can be determined by solving the optimization problem
\begin{equation}
\min_n T(n).
\end{equation}
Taking into account the condition
\begin{equation}\label{eq:solvability}
b^2 \geq \frac{2}{n} p_\text{RS} N (T_C + n T_\text{CP}) = \frac{2}{n} p_\text{RS} N T_c + 2 p_\text{RS} N T_\text{CP}
\end{equation}
required for the existence of a real solution, the minimization can be done analytically, yielding
\begin{equation}
n_\text{opt} = \frac{T_C\left(2p_\text{RS}\ N\ T_\text{CP} + b \sqrt{2 p_\text{RS}\ N\ T_\text{CP}}\right)}{2T_\text{CP}\left( b^2-2 N\ p_\text{RS}\ T_\text{CP}\right)}.
\end{equation}
 In this model the only influence of lossy compression is given by the time to read/write checkpoints and to recover data structures. While a simple model for time to checkpoint is given, e.g., in~\cite{CalCapOlsSniGro2019}, here we just exemplarily show the influence in Figure~\ref{fig:checkpointing}, by comparing different write/read times for~checkpointing.
Reducing checkpoint size by lossy compression, thus reducing $T_\text{CP}, T_\text{DS}$ has a small but noticeable effect on the overall runtime. Note that this model neglects the impact of inexact checkpoints on the re-computation time, which might increase, e.g., due to iterative methods requiring additional steps to reduce the compression error. For iterative linear solvers this is done in~\cite{TaoDiLiaCheCap2018}; a thorough analysis for the example of parallel-in-time simulation with hybrid parareal methods can be done along the lines of~\cite{FischerGoetschelWs2018}.

\begin{table}[H]
\centering
\caption{Notation used for optimal checkpointing. Where applicable, units are given in brackets.}
\label{tab:checkpointing-notation}
\begin{tabular}{cccc} \toprule
$n$ & number of checkpoints & $T_\text{C}$ & time for actual computation [s]\\
$N$ & number of compute cores &  $\ \, T_\text{CP}$ & time to write/read a checkpoint [s]\\ 
$p_\text{RS}$ & probability of failure & $T_\text{DS}$ & time to recover data structures [s]\\
&per unit time and core [1/s]&  $T_\text{R}$ & recovery time $=T_\text{CP}+T_\text{DS}$ [s]\\
$N_\text{RS}$ & number of restarts & $T_\text{RS}$ & time for restart [s]  \\
$T$ & overall runtime (wall clock) [s]  & $b$& $:=1-T_\text{R}p_\text{RS}N$ \\\bottomrule
\end{tabular}
\end{table}

\subsection{Postprocessing and Archiving}\label{sec:archiving}

Storage and postprocessing of results from large-scale simulations requires techniques to efficiently handle the vast amount of data. Assuming that computation requires higher accuracy than needed for postprocessing (due to, e.g., error propagation and accumulation over time steps), lossy compression will be beneficial here as well. In the following, we briefly present examples from three application~areas.

\subsubsection{Crash Simulation} Simulation is a standard tool in the automotive industry, e.g., for the simulation of crash tests. For archiving data generated by the most commonly used crash simulation programs, the lossy compression code FEMzip\footnote{\url{https://www.sidact.com/femzip0.html}; 
FEMzip is a registered trademark of Fraunhofer Gesellschaft, Munich.}~\cite{Thole2004} achieves compression factors of 10--20~\cite{TeranTholeLorentz2007,MerMueTho2015}, depending on the prescribed error tolerance. More recently, correlations between different simulation results were exploited by using a predictive principal component analysis to further increase the compression factor, reporting an increase by a factor of $4.4$ for a set of 14 simulation results~\cite{MerMueTho2015}.

\subsubsection{Weather and Climate} Today, prediction of weather and climate is one major use of  supercomputing facilities, with a tremendous amount of data to be stored ({In 2017, ECMWF's data archive grew by about 233 TB per day}\footnote{\url{https://www.ecmwf.int/en/computing/our-facilities/data-handling-system}}~\cite{DueLeuBau2019}). 
Thus, using compression on the whole I/O system (main memory, communication, storage) can significantly affect performance~\cite{KuhKunLud2016}. Typically, ensemble simulations are used, allowing to exploit correlations. Three~different approaches are investigated in~\cite{DueLeuBau2019}. The most successful one takes forecast uncertainties into account, such that higher precision is provided for less uncertain variables. Naturally, applying data compression should not introduce artifacts, or change, e.g., statistics of the outputs of weather and climate models. The impact of lossy compression on several postprocessing tasks is investigated in~\cite{BakHamMicXuEtAl2016,PoppickNardiFeldmanBakerHammerling2018}, e.g., whether artifacts due to lossy compression can be distinguished from the inherent variability of climate and weather simulation data. Here, avoiding smoothing of the data due to compression via transform coding can be important, favouring quantize-then-transform methods or simpler truncation approaches like fpzip~\cite{LindstromIsenburg2006}.

\subsubsection{Computational Fluid Dynamics} Solution statistics of turbulent flow are used in~\cite{OteVinMarLauSch2018} to assess error tolerances for lossy compression. Data reduction is performed by spatial transform coding with the discrete Legendre transform~\cite{MarSchFis2016}, which matches the spectral discretization on quadrilateral grids. For turbulence statistics of turbulent flow through a pipe, Otero et al.~\cite{OteVinMarLauSch2018} reported a reduction of $98\%$ for an admissible $L^2$ error level of $1\%$, which is the order of the typical statistical uncertainty in the evaluation of turbulence quantities from direct numerical simulation data.

\section{Conclusions}

Data compression methods are valuable tools for accelerating PDE solvers, addressing larger problems, or archiving computed solutions. Due to floating-point data being compressed, only lossy compression can be expected to achieve reasonable compression factors; this matches perfectly with the fact that PDE solvers incur discretization and truncation errors. An important aspect is to model and predict the impact of quantization or decimation errors on the ultimate use of the computed data, in order to be able to achieve high compression factors while meeting the accuracy requirements. This requires construction and use of fast and accurate a posteriori error estimators complementing analytical estimates, and remains one of the challenging future research directions.

Meeting such application-dependent accuracy requirements calls for problem-specific approaches to compression, e.g., in the form of transforms where the quantization of transform coefficients directly corresponds to the compression error in relevant spatially weighted Sobolev norms. Discretization-specific approaches that are able to exploit the known structure of spatial data layout, e.g., Cartesian grids or adaptively refined mesh hierarchies, are also necessary for achieving high compression factors.

Utility and complexity of such methods are largely dictated by their position in the memory hierarchy. Sophisticated compression schemes are available and regularly used for reducing the required storage capacity when archiving solutions. On the other hand, accelerating PDE solvers by data compression is still in the active research phase, facing the challenge that computational overhead for compression can thwart performance gains due to reduced data transmission time. Thus, simpler compression schemes dominate, in particular when addressing the memory wall. Consequently, only a moderate, but nevertheless consistent, benefit of compression has been shown in the literature.

The broad spectrum of partially contradicting requirements faced by compression schemes in PDE solvers suggests that no single compression approach will be able to cover the need, and that specialized and focused methods will increasingly be developed---a conclusion also drawn in~\cite{CappelloDiLiLiangGokTaoYoonWuAlexeevChong2019}.

The trend of growing disparity between computing power and bandwidth, which could be observed during the last three decades and will persist for the foreseeable future of hardware development, means that data compression methods will only become more important over time. Thus, we can expect to see a growing need for data compression in PDE solvers in the coming years.

\section*{Acknowledgements}
This work has been partially supported by the German Ministry for Eduction and Research (BMBF) under project grant 01IH16005 (HighPerMeshes). We thank Florian Wende for implementing mixed-precision preconditioners, Alexander Kammeyer for implementation and testing of checkpoint/restart, and Thomas Steinke for many helpful discussions.

\bibliographystyle{elsarticle-num}
\bibliography{paper}

\end{document}